
\newcount\mgnf\newcount\tipi\newcount\tipoformule\newcount\greco

\tipi=2          
\tipoformule=0   


\global\newcount\numsec
\global\newcount\numfor
\global\newcount\numtheo
\global\advance\numtheo by 1

\def\senondefinito#1{\expandafter\ifx\csname#1\endcsname\relax}

\def\SIA #1,#2,#3 {\senondefinito{#1#2}%
\expandafter\xdef\csname #1#2\endcsname{#3}\else
\write16{???? ma #1,#2 e' gia' stato definito !!!!} \fi}

\def\etichetta(#1){(\veroparagrafo.\veraformula)%
\SIA e,#1,(\veroparagrafo.\veraformula) %
\global\advance\numfor by 1%
\write15{\string\FU (#1){\equ(#1)}}%
\write16{ EQ #1 ==> \equ(#1) }}

\def\letichetta(#1){\veroparagrafo.\verotheo
\SIA e,#1,{\veroparagrafo.\verotheo}
\global\advance\numtheo by 1
\write15{\string\FU (#1){\equ(#1)}}
\write16{ Sta \equ(#1) == #1 }}

\def\tetichetta(#1){\veroparagrafo.\veraformula 
\SIA e,#1,{(\veroparagrafo.\veraformula)}
\global\advance\numfor by 1
\write15{\string\FU (#1){\equ(#1)}}
\write16{ tag #1 ==> \equ(#1)}}

\def\FU(#1)#2{\SIA fu,#1,#2 }

\def\etichettaa(#1){(A\veroparagrafo.\veraformula)%
\SIA e,#1,(A\veroparagrafo.\veraformula) %
\global\advance\numfor by 1%
\write15{\string\FU (#1){\equ(#1)}}%
\write16{ EQ #1 ==> \equ(#1) }}

\def\BOZZA{
\def\alato(##1){%
 {\rlap{\kern-\hsize\kern-1.4truecm{$\scriptstyle##1$}}}}%
\def\aolado(##1){%
 {
{
 \rlap{\kern-1.4truecm{$\scriptstyle##1$}}}}}
}

\def\alato(#1){}
\def\aolado(#1){}

\def\veroparagrafo{\number\numsec}
\def\veraformula{\number\numfor}
\def\verotheo{\number\numtheo}

\def\Eq(#1){\eqno{\etichetta(#1)\alato(#1)}}
\def\eq(#1){\etichetta(#1)\alato(#1)}
\def\leq(#1){\leqno{\aolado(#1)\etichetta(#1)}}
\def\teq(#1){\tag{\aolado(#1)\tetichetta(#1)\alato(#1)}}
\def\Eqa(#1){\eqno{\etichettaa(#1)\alato(#1)}}
\def\eqa(#1){\etichettaa(#1)\alato(#1)}
\def\eqv(#1){\senondefinito{fu#1}$\clubsuit$#1
\write16{#1 non e' (ancora) definito}%
\else\csname fu#1\endcsname\fi}
\def\equ(#1){\senondefinito{e#1}\eqv(#1)\else\csname e#1\endcsname\fi}

\def\Lemma(#1){\aolado(#1){\hskip.5cm\bf Lemma \letichetta(#1)}}
\def\Theorem(#1){\aolado(#1){\hskip.5cm\bf Theorem \letichetta(#1)}}
\def\Proposition(#1){\aolado(#1){\hskip.5cm\bf Proposition \letichetta(#1)}}
\def\Corollary(#1){\aolado(#1){\hskip.5cm\bf Corollary \letichetta(#1)}}
\def\Remark(#1){\noindent\aolado(#1){\it\bf Remark \letichetta(#1):}}
\def\Definition(#1){\aolado(#1){\hskip.5cm\bf Definition 
\letichetta(#1)}}
\def\Example(#1){\aolado(#1) {\it \bf Example \letichetta(#1)}}

\def\include#1{
\openin13=#1.aux \ifeof13 \relax \else
\input #1.aux \closein13 \fi}

\openin14=\jobname.aux \ifeof14 \relax \else
\input \jobname.aux \closein14 \fi
\openout15=\jobname.aux






\newcount\pgn \pgn=1
\def\foglio{\number\numsec:\number\pgn
\global\advance\pgn by 1}
\def\foglioa{A\number\numsec:\number\pgn
\global\advance\pgn by 1}

\footline={\rlap{\hbox{\copy200}}\hss\tenrm\folio\hss}

\def\TIPIO{
\font\setterm=amr7 
\def \settepunti{\def\rm{\fam0\setterm}
\textfont0=\setterm   
\normalbaselineskip=9pt\normalbaselines\rm
}\let\nota=\settepunti}

\def\TIPITOT{
\font\twelverm=cmr12
\font\twelvei=cmmi12
\font\twelvesy=cmsy10 scaled\magstep1
\font\twelveex=cmex10 scaled\magstep1
\font\twelveit=cmti12
\font\twelvett=cmtt12
\font\twelvebf=cmbx12
\font\twelvesl=cmsl12
\font\ninerm=cmr9
\font\ninesy=cmsy9
\font\eightrm=cmr8
\font\eighti=cmmi8
\font\eightsy=cmsy8
\font\eightbf=cmbx8
\font\eighttt=cmtt8
\font\eightsl=cmsl8
\font\eightit=cmti8
\font\sixrm=cmr6
\font\sixbf=cmbx6
\font\sixi=cmmi6
\font\sixsy=cmsy6
\font\twelvetruecmr=cmr10 scaled\magstep1
\font\twelvetruecmsy=cmsy10 scaled\magstep1
\font\tentruecmr=cmr10
\font\tentruecmsy=cmsy10
\font\eighttruecmr=cmr8
\font\eighttruecmsy=cmsy8
\font\seventruecmr=cmr7
\font\seventruecmsy=cmsy7
\font\sixtruecmr=cmr6
\font\sixtruecmsy=cmsy6
\font\fivetruecmr=cmr5
\font\fivetruecmsy=cmsy5
\textfont\truecmr=\tentruecmr
\scriptfont\truecmr=\seventruecmr
\scriptscriptfont\truecmr=\fivetruecmr
\textfont\truecmsy=\tentruecmsy
\scriptfont\truecmsy=\seventruecmsy
\scriptscriptfont\truecmr=\fivetruecmr
\scriptscriptfont\truecmsy=\fivetruecmsy
\def \eightpoint{\def\rm{\fam0\eightrm}
\textfont0=\eightrm \scriptfont0=\sixrm \scriptscriptfont0=\fiverm
\textfont1=\eighti \scriptfont1=\sixi   \scriptscriptfont1=\fivei
\textfont2=\eightsy \scriptfont2=\sixsy   \scriptscriptfont2=\fivesy
\textfont3=\tenex \scriptfont3=\tenex   \scriptscriptfont3=\tenex
\textfont\itfam=\eightit  \def\it{\fam\itfam\eightit}%
\textfont\slfam=\eightsl  \def\sl{\fam\slfam\eightsl}%
\textfont\ttfam=\eighttt  \def\tt{\fam\ttfam\eighttt}%
\textfont\bffam=\eightbf  \scriptfont\bffam=\sixbf
\scriptscriptfont\bffam=\fivebf  \def\bf{\fam\bffam\eightbf}%
\tt \ttglue=.5em plus.25em minus.15em
\setbox\strutbox=\hbox{\vrule height7pt depth2pt width0pt}%
\normalbaselineskip=9pt
\let\sc=\sixrm  \let\big=\eightbig  \normalbaselines\rm
\textfont\truecmr=\eighttruecmr
\scriptfont\truecmr=\sixtruecmr
\scriptscriptfont\truecmr=\fivetruecmr
\textfont\truecmsy=\eighttruecmsy
\scriptfont\truecmsy=\sixtruecmsy
}\let\nota=\eightpoint}

\newfam\msbfam   
\newfam\truecmr  
\newfam\truecmsy 
\newskip\ttglue
\ifnum\tipi=0\TIPIO \else\ifnum\tipi=1 \TIPI\else \TIPITOT\fi\fi

\def\a{\alpha}
\def\b{\beta}
\def\d{\delta}
\def\e{\epsilon}
\def\eps{\varepsilon}

\def\g{\gamma}

\def\l{\lambda}
\def\m{\mu}

\def\p{\pi}
\def\r{\rho}
\def\s{\sigma}

\def\z{\zeta}
\def\o{\omega}
\def\D{\Delta}
\def\L{\Lambda}
\def\G{\Gamma}
\def\O{\Omega}
\def\S{\Sigma}

\def\h{\eta}

\def\Sum{\Sigma}

\def\C{\hskip.2em\hbox{\rm C\kern-0.5em{I}\hskip.7em}}
\def\E{{I\kern-.25em{E}}}
\def\N{{I\kern-.25em{N}}}
\def\M{{I\kern-.25em{M}}}
\def\P{\hskip.2em\hbox{\rm P\kern-0.9em{I}\hskip.7em}}
\def\R{{I\kern-.25em{R}}}
\def\Z{{Z\kern-.425em{Z}}}
\def\1{{1\kern-.25em\hbox{\rm I}}}

\def\qed{{ \vrule height5pt width5pt depth0pt}\par\bigskip}

\def\AA{{\cal A}}
\def\BB{{\cal B}}

\def\EE{{\cal E}}

\def\GG{{\cal G}}

\def\KK{{\cal K}}

\def\LL{{\cal L}}
\def\ZZ{{\cal Z}}

\def\XX{{\cal X}}
\def\QQ{{\cal Q}}

\def\chap #1#2{\line{\ch #1\hfill}\numsec=#2\numfor=1}



\newcount\foot
\foot=1
\def\note#1{\footnote{${}^{\number\foot}$}{\ftn #1}\advance\foot by 1}
\def\tag #1{\eqno{\hbox{\rm(#1)}}}

\def\frac#1#2{{#1\over #2}}
\def\sfrac#1#2{{\textstyle{#1\over #2}}}
\def\text#1{\quad{\hbox{#1}}\quad}

\let\ra=\rightarrow

\def\proof{{\noindent \pr Proof: }}

\def\thanks{\noindent{\bf Aknowledgements: }}
\font\pr=cmbxsl10


\font\ch=cmbx12
\font\ftn=cmr8

\font\it=cmti10
\font\bf=cmbx10

%
\catcode`\X=12
\catcode`\@=11
\def\n@wcount{\alloc@0\count\countdef\insc@unt}
\def\n@wwrite{\alloc@7\write\chardef\sixt@@n}
\def\n@wread{\alloc@6\read\chardef\sixt@@n}
\def\crossrefs#1{\ifx\alltgs#1\let\tr@ce=\alltgs\else\def\tr@ce{#1,}\fi
   \n@wwrite\cit@tionsout\openout\cit@tionsout=\jobname.cit 
   \write\cit@tionsout{\tr@ce}\expandafter\setfl@gs\tr@ce,}
\def\setfl@gs#1,{\def\@{#1}\ifx\@\empty\let\next=\relax
   \else\let\next=\setfl@gs\expandafter\xdef
   \csname#1tr@cetrue\endcsname{}\fi\next}
\newcount\sectno\sectno=0\newcount\subsectno\subsectno=0\def\r@s@t{\relax}
\def\resetall{\global\advance\sectno by 1\subsectno=0
  \gdef\firstpart{\number\sectno}\r@s@t}
\def\resetsub{\global\advance\subsectno by 1
   \gdef\firstpart{\number\sectno.\number\subsectno}\r@s@t}
\def\v@idline{\par}\def\firstpart{\number\sectno}
\def\l@c@l#1X{\firstpart.#1}\def\gl@b@l#1X{#1}\def\t@d@l#1X{{}}
\def\m@ketag#1#2{\expandafter\n@wcount\csname#2tagno\endcsname
     \csname#2tagno\endcsname=0\let\tail=\alltgs\xdef\alltgs{\tail#2,}%
  \ifx#1\l@c@l\let\tail=\r@s@t\xdef\r@s@t{\csname#2tagno\endcsname=0\tail}\fi
   \expandafter\gdef\csname#2cite\endcsname##1{\expandafter
     \ifx\csname#2tag##1\endcsname\relax?\else{\rm\csname#2tag##1\endcsname}\fi
    \expandafter\ifx\csname#2tr@cetrue\endcsname\relax\else
     \write\cit@tionsout{#2tag ##1 cited on page \folio.}\fi}%
   \expandafter\gdef\csname#2page\endcsname##1{\expandafter
     \ifx\csname#2page##1\endcsname\relax?\else\csname#2page##1\endcsname\fi
     \expandafter\ifx\csname#2tr@cetrue\endcsname\relax\else
     \write\cit@tionsout{#2tag ##1 cited on page \folio.}\fi}%
   \expandafter\gdef\csname#2tag\endcsname##1{\global\advance
     \csname#2tagno\endcsname by 1%
   \expandafter\ifx\csname#2check##1\endcsname\relax\else%
\fi
   \expandafter\xdef\csname#2check##1\endcsname{}%
   \expandafter\xdef\csname#2tag##1\endcsname
     {#1\number\csname#2tagno\endcsnameX}%
   \write\t@gsout{#2tag ##1 assigned number \csname#2tag##1\endcsname\space
      on page \number\count0.}%
   \csname#2tag##1\endcsname}}%
\def\m@kecs #1tag #2 assigned number #3 on page #4.%
   {\expandafter\gdef\csname#1tag#2\endcsname{#3}
   \expandafter\gdef\csname#1page#2\endcsname{#4}}
\def\re@der{\ifeof\t@gsin\let\next=\relax\else
    \read\t@gsin to\t@gline\ifx\t@gline\v@idline\else
    \expandafter\m@kecs \t@gline\fi\let \next=\re@der\fi\next}
\def\t@gs#1{\def\alltgs{}\m@ketag#1e\m@ketag#1s\m@ketag\t@d@l p
    \m@ketag\gl@b@l r \n@wread\t@gsin\openin\t@gsin=\jobname.tgs \re@der
    \closein\t@gsin\n@wwrite\t@gsout\openout\t@gsout=\jobname.tgs }
\outer\def\localtags{\t@gs\l@c@l}
\outer\def\globaltags{\t@gs\gl@b@l}
\outer\def\newlocaltag#1{\m@ketag\l@c@l{#1}}
\outer\def\newglobaltag#1{\m@ketag\gl@b@l{#1}}

\def\t@gsoff#1,{\def\@{#1}\ifx\@\empty\let\next=\relax\else\let\next=\t@gsoff
   \expandafter\gdef\csname#1cite\endcsname{\relax}
   \expandafter\gdef\csname#1page\endcsname##1{?}
   \expandafter\gdef\csname#1tag\endcsname{\relax}\fi\next}
\def\verbatimtags{\let\ift@gs=\iffalse\ifx\alltgs\relax\else
   \expandafter\t@gsoff\alltgs,\fi}
\localtags
%
\global\newcount\numpunt
\magnification=\magstep1
\hoffset=0.cm
\baselineskip=14pt  
\parindent=00pt
\lineskip=4pt\lineskiplimit=0.1pt
\parskip=0.1pt plus1pt

\hyphenation{small}



\centerline {\bf Convergence to equilibrium for finite Markov processes,}
\centerline {\bf with application to the Random Energy Model.}
\vskip1cm
\centerline{ 
Pierre  Mathieu \footnote{$^1$}{\eightrm  CMI, Universit\'e de Provence,
39 Rue F. Joliot Curie, 13453 Marseille Cedex 13, France.\hfill\break
pmathieu$\!\,@$gyptis.univ-mrs.fr} and 
Pierre Picco \footnote{$^2$}{\eightrm CPT. CNRS Luminy, case 907, 
13288 Marseille Cedex 9, France.
Picco$\!\,@$cpt.univ-mrs.fr}.
}
\vskip.5truecm
\centerline{ CMI and CPT--CNRS, Marseille} 

\vskip.5cm


\vskip.5cm
{\bf Abstract:} {\it 
We estimate the distance in total variation between the law of 
a finite state Markov process at time $t$, 
starting from a given initial measure, and 
its unique invariant measure. 
We derive upper bounds for the time to reach the equilibrium. 
As an example of application we consider a special case of 
finite state Markov process in random environment: the Metropolis dynamics of 
the Random Energy Model.
We also study the process of the environment as seen from the process.
}

\vskip1truecm

{\bf Key Words:} {\it Finite Markov processes, Spectral gap, Poincar\'e inequality, 
Random Spin Systems, Spin glasses, Metropolis Dynamics.}
\vskip1truecm

{\bf AMS Classification Numbers}:  60K35, 82B44, 82D30, 82C44 
\vskip1truecm
{\bf Abbreviated title: } Finite Markov processes.
\vfill\eject


\def\XX{{\cal X}}

\chap{I. Introduction}1
\vskip.5truecm
\numsec=1
\numfor=1
\numtheo=1

We are interested in estimating the speed of convergence 
towards equilibrium for a finite and reversible Markov chain, a well studied 
problem in the theory of Markov chains, see [\rcite{SC}] for instance. 
Most, if not all results in this direction yield bounds on the distance 
to equilibrium which are uniform with respect to the initial distribution of the 
chain. In this paper, we shall rather derive estimates on mixing times that take 
into account the dependence on the initial law. 
As an example of application of our method, we  study the Metropolis dynamics of 
Derrida's Random Energy Model (REM.). 

Convergence times for the Metropolis dynamics of spin glasses were  
considered in [\rcite{M2}]. 
Let us note that the present paper was done simultaneously with  [\rcite{M2}] 
and  quoted therein as  [11] with a slightly different title. 
In [\rcite{M2}], estimates on the convergence time 
that depend on the initial law are given for models of spin glasses
such as the REM or the Sherrington-
Kirkpatrick model at high temperature. Three dynamics are considered: the random hoping 
time dynamics (RHT), 
the Glauber dynamics and the Metropolis dynamics. The initial configuration of the dynamics 
is always assumed to be chosen uniformly among all configurations. 

To compare the results obtained in the two articles, let us mention that 
 the starting points of the present article and [\rcite{M2}] are the same: the 
generalized Poincar\'e
inequalities that were introduced in [\rcite{M}], see section II here and in [\rcite{M2}]. 
However the way to estimate the associated constant $\LL_{\eta}(p)$, see \eqv(2.4) here and 
(2.2) there,  are completely different. We will come back to this point later. 

Since two  slightly 
different notions of convergence time  are used here and there, we first note that
in [\rcite{M2}], the time 
called $T^{\o}(c)$, is defined as in \eqv(P.1), with $c$ playing the r\^ole of $\e$.  
$T^{\o}(c)$ a priori depends on the realizations of 
the energies as the ${\o}$ emphasizes. In any case, the initial law, $\eta$,  is uniform.  

Here, for the Metropolis dynamics of the REM, the results  are given in term of a time 
denoted $T_N(\e,c,\eta)$ which is 
independent of the realizations of ${\o}$,  see \eqv(4.8).
It follows from the definition \eqv(4.8) that 
on a subset $\O_N$ of realizations of energies that has a  
probability larger than $1-e^{-cN}$ we have
$$
T_N(\e,c,\eta) \ge T^{\o}(\e)
$$
in particular this implies that, almost surely 
$$
\limsup_{N\rightarrow \infty} \frac 1{N}\log T^{\o}(\e) \le 
\limsup_{N\rightarrow \infty} \frac 1{N}\log T_N(\e,c,\eta) 
\Eq(P01)$$

We now recall some results from  [\rcite{FIKP}] and [\rcite{M2}] for the convergence time 
of the Metropolis dynamics of the REM. 
In [\rcite{M2}], it was proven that for
$\eta$ the uniform measure on $\{-1,+1\}^N$, we have, 
for almost all 
$\o$ 
$$
\limsup_{N\rightarrow \infty} \frac 1{N}\log T^{\o}(\e) \le 2\b^2 \,\,{\rm when}\,\, \b\le \b_c
\Eq(P02)
$$
and
$$
\limsup_{N\rightarrow \infty} \frac 1{N}\log T^{\o}(\e) \le 2\b\b_c \,\,{\rm when}\,\, \b\geq \b_c
\Eq(P03)
$$
(Remember that the free energy  and the mean energy per site converge  
for almost all $\o$ as it follows  from [\rcite{OP}]).  
Note however that using the  spectral gap estimates  for the Metropolis Dynamics of the REM
given in [\rcite{FIKP}], we  immediately get that,   for all $\b>0$, almost surely in $\o$ 
$$
\limsup_{N\rightarrow \infty} \frac 1{N}\log T^{\o}(\e) \le \b\b_c
\Eq(P04)
$$
and by checking all the probability estimates in [\rcite{FIKP}], we also have 
for all $\b>0$, for all $c>0$
$$
\limsup_{N\rightarrow \infty} \frac 1{N}\log T_N(\e,c,\eta) \le \b\b_c
\Eq(P041)
$$
Therefore, \eqv(P02) gives an better estimate  than \eqv(P04) only for $\b\le \b_c/2$.
For $\b>\b_c/2$, \eqv(P04) gives an  better estimate  than \eqv(P02) and \eqv(P03). 

Here we prove that for the Metropolis Dynamics of the REM, for all $\b\le \b_c$,
$$
\limsup_{N\rightarrow \infty} \frac 1{N}\log T_N(\e,c,\eta) \le \b^2
\Eq(P05)
$$
which together with \eqv(P01) and \eqv(P041)  gives for all $\b>0$ a better estimate  
than \eqv(P02) and \eqv(P03). 
Thus we have improved the results of [\rcite{M2}] in two ways: 
first we are using a more precise definition for the convergence time, second  
we won a factor $2$ in the upper bound for $\b\le \b_c$.

Note however that to get \eqv(P04) or \eqv(P041) a very careful analysis of 
optimization problems for paths on the weighted graph structure induced by 
the transition matrix of the dynamics was used. 
To prove \eqv(P05), a similar analysis is needed.  
Thus, using the specific paths constructed in [\rcite{FIKP}], instead of 
techniques based on estimates of the partition function as in [\rcite{M2}], 
leads, for the Metropolis 
dynamics of the REM, to an improvement  by a factor 2 in the estimates. 

We believe that the bound \eqv(P05) is sharp {\it i.e}
$\lim_{n\uparrow \infty}\frac 1{N}\log T_N(\e,c,\eta)= \b^2$ for $\b\le \b_c$.

We also believe that a similar analysis could be carried over for the Glauber dynamics, 
but the numerical factor in front of $\b^2$ in  \eqv(P05) would then be different. 
As far as the Random Hoping Time dynamics is concerned, it seems that the techniques 
of 
[\rcite{M2}]  directly lead to upper bounds of the correct order. Note however that 
the RHT dynamics has a much simpler structure than the Metropolis or Glauber ones. Indeed 
the RHT dynamics is nothing but a time-changed standart random walk on configuration space. 
The sequence of the different states visited by the process is independent of the 
Hamiltonian.  
On one hand, this feature very much simplifies the geometry. On the other hand, 
physicists believe that 
the evolution of the process should rather look like  a random perturbation of  
the steepest gradient dynamical system. The RHT dynamics displays un-physical features.

The organization of the paper is as follows: part II and II deal with general 
reversible Markov chains on a finite set.  In part II, 
we define generalized Poincar\'e inequalities and show that they
control the decay of the semi-group (Theorem \eqv(theo1)). Then we derive geometric
estimates for the generalized Poincar\'e constants (Theorem \eqv(theo2)). Part III
contains an application of these results in a case where the state space can be
splitted into two components: 'good' and 'bad' points. For the reader's convenience, 
we decided to give self-contained proofs of our results at the risk of repeating 
arguments already used in [\rcite{M}], [\rcite{M1}] or [\rcite{M2}] .  

Although we shall not directly
use the results of part III to study the R.E.M., the strategy will be the same. Only
technical aspects make the computation for the R.E.M. a little longer than the proof in
part III. In part IV, we precisely define the R.E.M. and state our bounds for the
thermalization time ( Theorems \eqv(theo5) and \eqv(theo6)). Then we proceed to the
proofs. 
In part V, we extend our results to the process of the environment as seen from
the particle. This section is similar to the section 3 of [\rcite{M2}] 
with more pedagogical details  on the construction of the process.
We then show that the equilibrium time also satisfies \eqv(1.22). Part
VI contains the proof of some static estimates on the R.E.M. that we needed in the
previous parts.


\vskip1cm

\chap{II. Generalized Poincar\'e inequalities}2
\vskip.5truecm
\numsec=2
\numfor=1
\numtheo=1

Let $X=(X_t)_{t\ge 0}$ be an homogeneous Markov 
process on a finite state space, $\XX$. 
We assume that there is a unique invariant, ergodic probability measure 
for $X$, say $\p$. We further assume that 
$\p$ charges every point in $\XX$ and that it is reversible. 
Let $\h$ be some probability measure on $\XX$ and 
call $\LL_\h(X_t)$ the law of $X_t$ when the initial law 
is $\h$. We wish to bound $d_{TV}(\LL_\h(X_t), \p)$, 
the distance in total variation between the law of $X$ 
at time $t$ and the equilibrium law $\p$. More precisely, 
we would like to obtain an upper bound in terms 
of the geometry of the Markov process $X$ i.e. in terms 
of the geometry of the graph structure 
induced by the transition matrix 
on the state space.

It is well known that one can use Poincar\'e inequalities 
to bound $d_{VT}(\LL_\h(X_t), \p)$. Indeed 
calling $\l$ the spectral gap of the generator of $X$ 
(which is a symmetric matrix since we have assumed 
that $\p$ is reversible), we have, for any real valued function $f$ 
defined on $\XX$ and for any $t\ge 0$, 
$$\p[(P_tf-\p(f))^2]\le e^{-2\l t} \p(f^2) \Eq(1.1)$$
where $P_t$ denotes the semi-group  i.e. $P_tf(x)=E_x[f(X_t)]$. 
From \eqv(1.1), it immediately follows that 
$$
max _{x\in\XX}d_{TV}(\LL_x(X_t), \p)
\le \sqrt{\frac1{\p_*}} e^{-\l t} \Eq(1.2)$$
where $\LL_x(X_t)$ is the law of $X_t$ when the initial law 
is a Dirac mass at the point $x\in\XX$ and 
$\p_*=\min_{x\in\XX} \p(x)$. 
It now remains to estimate $\l$ in terms of the geometry of $X$. 
Such bounds exist, they rely on 
Poincar\'e inequalities: assume that for some constant $a>0$ and 
any function $f$ with $\p(f)=0$, we have:
$$\p(f^2)\le a \EE(f,f)\Eq(1.21)$$
then $1/\l\le a$. Here $\EE$ is  the Dirichlet form of $X$. 
From \eqv(1.21) one can deduce lower bounds of $\l$ in terms of 
optimization problems for paths on the weighted graph structure 
induced by the transition matrix of $X$ on $\XX$ (See [\rcite{SC}] 
and the references therein).  
\eqv(1.2) might be sharp or not depending on $X$. Many efforts were 
recently made to improve 
\eqv(1.2). More precise bounds can be obtained replacing the 
Poincar\'e inequality by 
 more sophisticated functional inequalities 
such as Log-Sobolev, Sobolev or Nash inequalities. 
We refer to [\rcite{SC}] for a detailed discussion 
of this topic. In all  cases, one estimates  
$max _{x\in\XX}d_{TV}(\LL_x(X_t), \p)$ i.e. the speed of 
convergence to equilibrium starting from the worst initial point.

We look for estimates of 
$d_{TV}(\LL_\h(X_t), \p)$ that should depend on $\h$. 
This paper is an  attempt to adapt the strategy of the Poincar\'e 
inequality in this context: for each initial law $\h$,
 we introduce a family of functional inequalities, quite 
similar to the Poincar\'e one, and prove that 
they allow one to control the distance to equilibrium. 
We call these inequalities {\it generalized Poincar\'e inequalities}. 
We then derive geometric bounds for the 
constants involved in these inequalities in the spirit of [\rcite{SC}].

Let 
$(K(x,y), (x,y)\in\XX\hbox{x}\XX)$ be the transition matrix of the   
Markov process $X$. 
Since we assume that the 
measure $\p$ is reversible, the kernel $k(x,y)=K(x,y)/\p(x)$ is symmetric, 
i.e. $k(x,y)=k(y,x)$. 
Let $P_tf(x)=E_x[f(X_t)]$ denote the semi-group associated to $X$. 

For functions $f$ and $g$ defined on $\XX$, let 
$$\EE(f,g)==\frac 12 \sum_{x,y}(f(x)-f(y))(g(x)-g(y)) k(x,y)\p(x)\p(y)\Eq(2.1)$$
be the Dirichlet form of $X$. 
For any edge $e=(x,y)\in\XX\hbox{x}\XX$, 
let $Q(e)=k(x,y)\p(x)\p(y)$. Also define $d_ef=f(x)-f(y)$. Then 
\eqv(2.1) can be re-written as 
$$\EE(f,g)=\frac 12 \sum_{e\in\XX^2} Q(e) d_ef d_eg \Eq(2.2)$$

For $p\in]0,1]$, let us define the following constants:
$$\LL(p)=\inf_{f\ s.t.\ \p(f)=0} 
\frac{\EE(f,f)\Vert f\Vert_\infty^{(2-2p)/p}}{\p(\vert f\vert)^{2/p}} \Eq(2.3)$$
and, for a probability measure on $\XX$, say $\h$, 
$$\LL_\h (p)=\inf_{f\ s.t.\ \p(f)=0} 
\frac{\EE(f,f)\Vert f\Vert_\infty^{(2-2p)/p}}{\h(\vert f\vert)^{2/p}} \Eq(2.4)$$
Clearly $\LL(p)=\LL_\p(p)$. 
H\"older's inequality implies that the function $p\ra\LL_\h(p)$ is decreasing and 
that $\LL(p)\ge \l$ for any $p$. 
(Remember that $\l$ denotes the spectral gap of the generator of $X$). 

To measure the time it takes for the process 
to reach equilibrium, we define the following quantities: 
$$d_\eta(t)=\sup_{s\ge t}\sup_{f;\Vert f\Vert_\infty\le 1} 
\eta(\vert P_sf-\pi(f)\vert)$$
and, for any $\eps>0$, 
$$T_\h(\eps)=\inf\{t>0\ s.t.\  d_\eta(t) \le\eps\}\Eq(P.1)$$
Note that $d_{TV}(\LL_\h(X_s),\p)\le d_\eta(t)$ for all $s\ge t$.

{\bf Remark}: let 
$$\L(p)=\inf_{f\ s.t.\ \p(f)=0} \frac{\EE(f,f)}{\p(\vert
f\vert^p)^{2/p}}$$ 
Then, as a consequence of H\"older's inequality, $\L(p)\le\LL(p)$ for $p\in]0,1]$. 
Also $\L(1)=\LL(1)$. The constants $\L(p)$ and $\LL(p)$ have 
already been introduced in [\rcite{M}]. 
(In the notation of [\rcite{M}], $\LL(p)$ is denoted 
$\KK(p/(1-p),+\infty)$). 
It follows from the results of [\rcite{M}], 
that $\LL(p)$ can also be defined in terms 
of the  capacity associated to $\EE$ 
and different estimates of hitting times can be derived 
in terms of $\LL(p)$.\footnote{$^*$}
{\eightrm We take this opportunity to warm the reader that the results of 
part II in [\rcite{M}] are false.} 
We also have $\L(2)=\l$ and $\L(p)\ge \l$ for any $p\in]0,1]$. 
Because of the similarity of the 
definition of $\L(p)$ and the Poincar\'e inequality, 
we call the inequality $\L(p)\ge a$ for some $a>0$, 
a ''generalized Poincar\'e inequality'', 
although there is no spectral interpretation.

{\bf \Theorem(theo1)}: {\it let $p\in]0,1]$ and $p'\in]0,1]$. There exists a universal 
function of $(p,p')$, $C_{p,p'}$, such that, for any 
probability measure $\h$ and any $t>0$, 
$$d_{\eta}(t)\le C_{p,p'} \LL_\h(p')^{-p'/2} 
\LL(p)^{-pp'/(4-2p)}t^{-p'/(2-p)} \Eq(2.5)$$
$C_{p,p'}=e^{-p'/2} (p/(2-p))^{pp'/(4-2p)}$ would do.
As a consequence, for any $\eps>0$, we have
$$T_\h(\eps)\le {\tilde C}_p \LL_\h(p')^{-(2-p)/2}\LL(p)^{-p/2} 
\eps^{-(2-p)/p'}\Eq(2.51)$$
where ${\tilde C}_p= e^{-(2-p)/2} (p/(2-p))^{p/2}$.
}

{\proof}: we shall prove that, for any function $f$ with $\p(f)=0$, then 
$$\h(\vert P_t f\vert)\le
 C_{p,p'} \LL_\h(p')^{-p'/2} \LL(p)^{-pp'/(4-2p)}
t^{-p'/(2-p)}\Vert f\Vert_\infty \Eq(2.6)$$
with $C_{p,p'}=e^{-p'/2}  (p/(2-p))^{pp'/(4-2p)}$. 
\eqv(2.6) implies \eqv(2.5).

{\it Step 1}: define 
$$\KK(p)=\inf_{f\ s.t.\ \p(f)=0} 
\frac{\EE(f,f)\Vert f\Vert_\infty ^{(4-2p)/p}}
{\p(\vert f\vert^2)^{2/p}} \Eq(2.7)$$
From H\"older's inequality, we deduce that $\KK(p)\ge \LL(p)$. 
Let $f$ be s.t. $\p(f)=0$. We claim that 
$$\p[(P_t f)^2]\le  
(\frac {4-2p}p)^{-p/(2-p)}(\KK(p)t)^{-p/(2-p)}\Vert f\Vert_\infty^2\Eq(2.8)$$
Then \eqv(2.8) will also hold with $\KK(p)$ replaced by $\LL(p)$.

{\it Proof of \eqv(2.8)}: let $f(t)=\p[(P_t f)^2]$. 
Then $f'(t)=-2\EE(P_tf,P_tf)$. 
By definition of $\KK(p)$, we have:
$$f'(t)\le -2\KK(p) \frac{f(t)^{2/p}}{\Vert P_t f\Vert_\infty^{(4-2p)/p}}$$
Since $P_t$ is a contraction in $L_\infty$, we also have
$$f'(t)\le -2\KK(p) \frac{f(t)^{2/p}}{\Vert f\Vert_\infty^{(4-2p)/p}}$$
Integrating this last inequality, we get

$$\eqalign{
f(t)^{1-2/p}&\ge f(0)^{1-2/p}+ 
\frac{4-2p}p t \frac{\KK(p)}{\Vert f\Vert_\infty^{(4-2p)/p}}\cr
 &\ge  \frac{4-2p}p t \frac{\KK(p)}{\Vert f\Vert_\infty^{(4-2p)/p}}  }
$$
which implies \eqv(2.8).
\qed
{\it Step 2}: there exists a universal 
constant $C$ s.t. for any function $f$ and any $t>0$ we have 
$$\EE(P_t f,P_t f)\le (C/t) \p[f^2]\Eq(2.9)$$
($C=1/(2e)$ would do.) 

{\it Proof}: for all $\m\ge 0$ and $t>0$, 
we have  $\m e^{-2\m t}\le C/t$. Use this inequality and 
a spectral decomposition of the Dirichlet form $\EE$ to deduce \eqv(2.9).\qed

{\it Step 3}: we finish the proof of \eqv(2.6). 
By definition of $\LL_\h(p')$, we have:
$$\h(\vert P_t f\vert)^{2/p'}\le  
\LL_\h(p')^{-1} \EE(P_tf,P_tf)\Vert P_t f\Vert_\infty^{(2-2p')/p'}$$
Using \eqv(2.9), the semi-group 
property: $P_t=P_{t/2}P_{t/2}$, and 
the fact that $P_t$ is a contraction in $L_\infty$, we get that 
$$\h(\vert P_t f\vert)^{2/p'}\le C 
\frac 2{\LL_\h(p')t} 
\p(\vert P_{t/2} f\vert^2) \Vert f\Vert_\infty^{(2-2p')/p'}$$
Using \eqv(2.7) ( with $\LL(p)$ instead of $\KK(p)$ ), we get that 
$$\h(\vert P_t f\vert)^{2/p'}\le 
2C (\frac {4-2p}p)^{-p/(2-p)} 
(\LL_\h(p') t)^{-1}  (\frac 2{\LL(p)t})^{p/(2-p)} \Vert f\Vert_\infty^{2/p'}$$
\qed

{\bf Remarks}: 

(i) Depending on the concrete example under consideration, the sharpness of the 
bound \eqv(2.5) ranges from good to extremely bad. Let us just outline one example where 
Theorem \eqv(theo1) leads to a very bad estimate: we consider the usual random walk on the 
discrete cube $\XX=\{-1,+1\}^N$. Then $Q(e)=1/(N2^N)$, 
for any edge between two nearest neighbours in $\XX$. 
Choose for $\eta$ a Dirac mass, say $\eta=\d_a$. Using the test function 
$f=\d_a-\pi(a)$ in formula \eqv(2.4), we get that, for large enough $N$, 
$$\LL_\eta(p)\le 2^{2/p-N}$$
Therefore \eqv(2.51) would lead to the conclusion that the process reaches equilibrium 
in a time shorter than $\exp(cN)$, whereas the true value of $T_\eta(\eps)$ is known to be 
of order 
$ N\log N$. We will see with the R.E.M.  an example where Theorem \eqv(theo1) leads 
to more interesting conclusions.

There is one situation in which \eqv(2.5) is not so far from being sharp: 
assume that $\eta=\pi$. Let $a$ be such that, for any function $f$ with $\pi(f)=0$, and 
for any time $t>0$, we have 
$$\pi(\vert P_t f\vert]\le (\frac a t)^{\frac p{2-p}}\Vert f\Vert_\infty \Eq(r1)$$
By interpolation, \eqv(r1) implies that 
$$\pi[(P_t f)^2]\le (\frac a t)^{\frac p{2-p}} \Vert f\Vert_\infty^2 \Eq(r2)$$
Use now the inequality
$$\pi[f^2]-\pi[(P_t f)^2]=\int_0^t2\EE(P_sf,P_sf)ds\le2t\EE(P_tf,P_tf)$$
to get that 
$$\pi(f^2)\le 2t\EE(P_tf,P_tf)+(\frac a t)^{\frac p{2-p}} \Vert f\Vert_\infty^2$$
Choosing the best value for $t$, we obtain the inequality:
$$\pi(f^2)\le C_pa^{\frac p2} \Vert f\Vert_\infty^{2-p} (\EE(f,f))^{\frac p2}$$
,where $C_p$ is some universal function of $p$.
In other words we have proved that $1/\KK(p)\le C_p a$, i.e \eqv(2.8) is sharp, 
up to multiplicative constants.

(ii) We derive estimates of the eigenvectors of $\EE$ in terms of $\LL(p)$. Following 
the terminology of [\rcite{M}], let us define 
$$\KK_2(p)=\inf_{f\ s.t.\ \pi(f)=0} 
 \frac { \EE(f,f) \pi(f^2)^{(1-p)/2p} } {\pi(\vert f\vert)^{(p+1)/p} } \Eq(r3)$$ 
 
It follows from Proposition 1, Proposition 2 and Theorem 1 in [\rcite{M}] that, for any 
$p'<p$, there exists a constant $C_{p,p'}$ such that $\LL(p)\le C_{p,p'} \KK_2(p')$.
 
Let now $l$ be an eigenvalue of $\EE$ and $\phi$ be the corresponding eigenvector. 
We assume that $l\not=0$ ($\phi$ is not constant),  and $\pi[\phi^2]=1$. Using $f=\phi$ in \eqv(r3) and 
$\EE(\phi,\phi)=l$, we obtain that $\KK_2(p)\le l/\pi(\vert \phi\vert)^{(p+1)/p}$. 
Replacing $\KK_2(p)$ by $\LL(p)$, we therefore have:
$$\pi(\vert \phi\vert)\le C_{p,p'} (\frac l{\LL(p)})^{p'/(1+p')} \Eq(r4)$$
for any $p'<p$. 

\eqv(r4) implies that, if $l$ is much smaller than $\LL(p)$, then $\pi(\vert\phi\vert)$ is 
small i.e. the function $\phi$ is very concentrated on its support. Since $l\ge \l$, 
where 
$\l$ is the spectral gap, this situation can  occur only if, for some $p$, $\l<<\LL(p)$. 
This will be the case for Metropolis dynamics of the R.E.M. at high temperature and 
we shall use 
\eqv(r4)  to prove that the first 
eigenvector of the dynamics is degenerate.

{\bf Geometric estimates}: a path $\g$ in $\XX$ 
is a sequence of vertices $\g=(x_0,...,x_k)$. 
Equivalently, $\g$ can be viewed as a sequence of bounds 
$\g=(e_1, ...,e_k)$ with 
$e_i=(x_{i-1},x_i)$. The length of $\g$ is $\vert\g\vert=k$. For $x,y\in\XX$, let 
$\G(x,y)$ be the set of all paths 
$\g=(x_0,...,x_k)$ with $x_0=x$ and $x_k=y$ and $k(x_{i-1},x_i)\not=0$ for 
all $i=1...k$. For each $x\not=y\in\XX$, let us 
choose one path, say $\g(x,y)\in \G(x,y)$.  
Since we have assumed that $\p$ is ergodic and 
charges all points in $\XX$, $X$ is 
irreducible and therefore $\G(x,y)$ is always non empty.

{\bf \Theorem(theo2)}: {\it (i) Let $p\in]0,1[$. Let $\l(x)$ and $\m(x)$ be 
two positive functions on $\XX$. We have
$$\eqalign{
 \frac 1{\LL_\h(p)}\le 2^{2/p-1} 
&\left( \sum \p(x) \l(x)^{p/(1-p)} 
\sum \h(y) \m(y)^{p/(1-p)} \right)^{(2-2p)/p}\cr 
&\left( \sum_{e\ s.t.\ Q(e)\not=0} \frac 1{Q(e)} 
(\sum_{x,y\ s.t.\ e\in\g(x,y)} \frac{\p(x)\h(y)}{\l(x)\m(y)} )^2\right) 
}\Eq(2.20)$$
(ii) 
$$ \frac 1{\LL_\h(1)}\le  
2 \left( \sum_{e\ s.t.\ Q(e)\not=0} 
\frac 1{Q(e)} (\sum_{x,y\ s.t.\ e\in\g(x,y)} \p(x)\h(y))^2\right) 
\Eq(2.21)$$
}

{\it Comments}: let us recall from [\rcite{SC}] the following estimate of the 
spectral gap:
$$\frac 1\l\le 
\max_{e\ s.t.\ Q(e)\not=0} \frac 1{Q(e)}\sum_{x,y\ s.t.\ e\in\g(x,y)}
\vert \g(x,y)\vert \p(x)\p(y) \Eq(2.211)$$

{\proof}: (ii) follows from (i): choose $\l(x)=\m(x)=1$, and let $p$ tend to $1$.

Let $f$ be a function s.t. 
$\p(f)=0$. Note that $f(y)-f(x)=\sum_{e\in\g(x,y)} d_ef$. Therefore  
$$\eqalign{
\Sum_y \h(y)\vert f(y)\vert = 
& \Sum_y \h(y)\vert f(y)\vert-\Sum_x \p(x) f(x)\cr
&\le \Sum_{x,y} \h(y)\p(x) \vert f(x)-f(y)\vert\cr
&=\Sum_{x,y} \h(y)\p(x)\vert\Sum_{e\in\g(x,y)}d_e f\vert\cr
&\le 2^{1-p}\Sum_{x,y} \h(y)\p(x) \vert 
\Sum_{e\in\g(x,y)} d_e f\vert ^p \Vert f\Vert _\infty^{1-p}\cr
&= 2^{1-p} \Vert f\Vert _\infty^{1-p} 
\Sum_{x,y} \h(y)\p(x) [\l(x)\m(y)\vert\Sum_{e\in\g(x,y)} d_e f\vert]^p
\l(x)^{-p}\l(y)^{-p}
}$$

We apply H\"older's inequality to get 
$$\eqalign{
&\Sum_y \h(y)\vert f(y)\vert\cr
&\le  2^{1-p} \Vert f\Vert _\infty^{1-p} 
\left( \Sum_{x,y}\p(x)\h(y) \l(x)^{p/(1-p)} \m(y)^{p/(1-p)} \right)^{1-p}
 \left(\Sum_{x,y} \frac{\p(x)\h(y)}{\l(x)\m(y)} 
\vert\Sum_e d_ef\vert\right)^p \cr
&\le 2^{1-p} \Vert f\Vert _\infty^{1-p} 
\left( \Sum_{x,y}\p(x)\h(y) 
\l(x)^{p/(1-p)} \m(y)^{p/(1-p)} \right)^{1-p} \times \cr
&\quad\quad \times  \left(\Sum_e \vert d_e f\vert 
\Sum_{x,y\ s.t.\ e\in\g(x,y)} \frac{\p(x)\h(y)}{\l(x)\m(y)} \right)^p 
}$$
Applying once more H\"older's inequality, we get 
$$\eqalign{
&\Sum_y \h(y)\vert f(y)\vert\cr
&\le2^{1-p} \Vert f\Vert _\infty^{1-p} 
\left( \Sum_{x,y}\p(x)\h(y) \l(x)^{p/(1-p)} \m(y)^{p/(1-p)} \right)^{1-p}\times\cr
&\quad\quad \times \left( \Sum_e (d_ef)^2 Q(e)\right)^{p/2} 
\left( \Sum_e \frac1{Q(e)} 
(\Sum_{x,y\ s.t.\ e\in\g(x,y)} \frac{\p(x)\h(y)}{\l(x)\m(y)} )^2\right)^{p/2}
}\Eq(2.22)$$
Replacing $\Sum_e Q(e) 
\vert d_e f\vert^2$ by $2\EE(f,f)$, we get the desired result.
\qed


\vskip1cm

\chap{III. Applications}3
\vskip.5truecm
\numsec=3
\numfor=1
\numtheo=1

This part of the paper mainly has 
a pedagogical aim. We shall illustrate how one can 
use the results of part II 
in a concrete situation. An even more concrete example of application 
will be given in the next part with the R.E.M. 

Comparing \eqv(2.21) and \eqv(2.211), one sees that the gain in 
using generalized spectral gap inequalities 
instead of the usual spectral gap inequality 
is that we can now afford having 
some ''very bad sites'' since we replaced a ''max'' over edges $e$ 
by a sum. Besides formula \eqv(2.20) 
gives us the possibility of 'killing' these bad points 
by choosing $\l$ and $\m$. 
To illustrate the way it works, let us 
assume that the state space $\XX$ can be divided 
into two disjoint sets, $B$ and $G$. 
'B' stands for 'bad'. Points in $B$ are supposed to be pathological 
and we do not expect them to play any role on the speed of convergence when 
the initial measure is smooth enough. 

The next Theorem states a lower bound for $\LL_\h(p)$ which is valid for any partition
of $\XX$ into two sets 
$B$ and $G$, but \eqv(3.1) is useful only if,  
firstly, we assume that the measure of $B$ is small both for 
$\p$ and $\h$ and besides we also assume somehow that the hitting time of 
$B$ is large i.e. the 
weights $Q(e)$ for those edges $e$ that touch $B$ are not too small.

Let us introduce some notation: 
$$\g^*=\sup_{x,y\in\XX}\vert \g(x,y)\vert$$
$$\BB=\{e\in \XX\hbox{x}\XX\ s.t.\ there\ exist\ 
x\ and\ y\ s.t.\ e\in\g(x,y)\ and\ x\in B\ or\ y\in B\}$$
In $\BB$ are edges $e\in\g(x,y)$ with both $x$ and $y$ in $B$. 

{\bf \Theorem(theo3)}: {\it for any $p\in]0,1]$, for any probability measure $\h$
$$\eqalign{
\frac 1{\LL_\h(p)}\le
2^{6/p-3}  \{&
\g^*\sup_{e\ s.t.\ Q(e)\not=0} 
\left(\frac 1 {Q(e)} 
\sum_{x\in G, y\in G\ s.t.\ e\in\g(x,y)} \p(x)\h(y)\right)\cr
&+ 2\left(\sum_{e\in \BB} 
\frac 1 {Q(e)}\right) \left( \p(B)^{2/p}+\h(B)^{2/p}\right)
\}}  \Eq(3.1)$$
}

{\proof}: let $p\in]0,1[$. The proof for $p=1$ is simpler and we 
leave it to the reader. 
Let us choose $\l$ and $\m$ as follows: $\l(x)=\m(x)=1$ for $x\in G$, 
$\l(x)=\p(B)^{1-1/p}$ for $x\in B$ and 
$\m(x)= \h(B)^{1-1/p}$ for $x\in B$. 
Then 
$$\eqalign{
&\sum \p(x)\l(x)^{p/(1-p)}\cr
=&\p(G)+1\le 2
}$$
The same holds for $\sum\h(y)\m(y)^{1/(1-p)}$. Therefore 
$$\frac 1 {\LL_\h(p)}
\le 2^{6/p-5} \sum_e \frac 1 {Q(e)} 
(\sum_{x,y\ s.t.\ e\in\g(x,y)} \frac{\p(x)\h(y)}{\l(x)\m(y)})^2
\Eq(3.2)$$

We compute the sum in \eqv(3.2) considering separately the
 cases $(x,y)\in G\hbox{x}G$, $(x,y)\in B\hbox{x}B$, 
$(x,y)\in G\hbox{x}B$ and $(x,y)\in B\hbox{x}G$. 
Since $\l=\m=1$ on $G$, the first term is bounded by 
$$\eqalign{
&\sum_e \frac 1{Q(e)} (\sum_{x,y\in G\ s.t.\ e\in\g(x,y)} \p(x)\h(y))^2\cr
\le &\left(\sup_e \frac 1{Q(e)} 
\sum_{x,y\in G\ s.t.\ e\in\g(x,y)} \p(x)\h(y)\right)
\left(\sum_e\sum_{x,y\in\XX\ s.t.\ e\in\g(x,y)} \p(x)\h(y)\right)\cr
=&\left(\sup_e \frac 1{Q(e)} \sum_{x,y\in G\ s.t.\ e\in\g(x,y)} \p(x)\h(y)\right)
\left(\sum_{x,y}\vert \g(x,y)\vert  \p(x)\h(y)\right)\cr
\le & \left(\sup_e \frac 1{Q(e)} 
\sum_{x,y\in G\ s.t.\ e\in\g(x,y)} \p(x)\h(y)\right)
\g^*}
\Eq(3.3)$$

The term corresponding to the case $(x,y)\in B\hbox{x}B$ is bounded by 
$$\eqalign{
& \sum_{e\in\BB} \frac 1{Q(e)} 
(\sum_{x,y\in B} \p(x) \h(y) \p(B)^{1/p-1}\h(B)^{1/p-1})^2\cr
\le & \left( \sum_{e\in\BB} \frac 1{Q(e)} \right) \p(B)^{2/p}\h(B)^{2/p}\cr
\le & \left( \sum_{e\in\BB} \frac 1{Q(e)} \right) ( \p(B)^{2/p}+\h(B)^{2/p})
}\Eq(3.4)$$

The term corresponding to the case $(x,y)\in G\hbox{x}B$ is bounded by 
$$\eqalign{
& \sum_{e\in\BB} \frac 1{Q(e)} \p(G)^2 \h(B)^{2/p} \cr
\le &  \sum_{e\in\BB} \frac 1{Q(e)} \h(B)^{2/p}
}\Eq(3.5)$$
Similarly the contribution of  $(x,y)\in B\hbox{x}G$ is bounded by 
$$\sum_{e\in\BB} \frac 1{Q(e)} \p(B)^{2/p}\Eq(3.6)$$
Inserting these bounds 
in \eqv(3.2) leads to the statement of Theorem \eqv(theo3).\qed

In the preceding Theorem, we chose the same 'bad' set for 
both measures $\p$ and $\h$. We now describe a 
slightly more sophisticated version of Theorem \eqv(theo3) 
obtained when choosing a different bad set for $\p$ 
and $\h$. Let us therefore assume 
that $\XX$ can be split into the disjoint union of two 
sets $B_\h$ and $G_\h$. $B_\h$ might differ from $B$. We modify 
the definition of $\BB$ accordingly: 
$$\BB=\{e\in \XX\hbox{x}\XX\ s.t.\ there\ exist\ 
x\ and\ y\ s.t.\ e\in\g(x,y)\ and\ x\in B\ or\ y\in B_\h\}$$
The proof of the following claim is identical to the proof of Theorem 
\eqv(theo3):

{\bf \Theorem(theo4)}: {\it for any $p\in]0,1]$, 
for any probability measure $\h$, 
any partitions $\XX=B\cup G=B_\h\cup G_\h$, 
we have
$$\eqalign{
\frac 1{\LL_\h(p)}\le
2^{6/p-3}  \{&
\g^*\sup_{e\ s.t.\ Q(e)\not=0} \left(\frac 1 {Q(e)} 
\sum_{x\in G, y\in G_\h\ s.t.\ e\in\g(x,y)} \p(x)\h(y)\right)\cr
&+ 2\left(\sum_{e\in \BB} \frac 1 {Q(e)}\right) 
\left( \p(B)^{2/p}+\h(B_\h)^{2/p}\right)
\}}  \Eq(3.7)$$
}

{\proof}: choose $\l(x)=1$ for $x\in G$, $\m(x)=1$ for $x\in G_\h$ and 
$\l(x)=\p(B)^{1-1/p}$ for $x\in B$, $\m(x)=\h(B_\h)^{1-1/p}$ for $x\in B_\h$. 
Then proceed as in the proof 
of Theorem \eqv(theo3).\qed

Finally let us mention that even more elaborated bounds can be obtained: 
we could distinguish bounds in $\BB$ 
linking sites $(x,y)$ with $(x,y)\in G\hbox{x} B_\h$,  $(x,y)\in B\hbox{x}G$ and 
$(x,y)\in B\hbox{x} B_\h$. We could also introduce 'weights' on bounds.
 We could choose a 'flow' 
of paths rather that picking a single path from $x$ to $y$. 
If necessary, one can also use these three 
tricks at the same time. We refer to Chapter 3 in [\rcite{SC}] 
for the notions of 'weights' and 'flow' or 
even 'generalized weights'.


\def\XX{{\cal X}}
\vskip1cm
\chap{IV. Dynamical phase transition for the  REM}4
\vskip.5truecm
\numsec=4
\numfor=1
\numtheo=1 

Before stating our result, let us recall the definition and some known facts 
on the R.E.M. 

{\bf Derrida's Random Energy Model}: 
The REM was introduced by Derrida 
[\rcite{B1},\rcite{B2}]
as the simplest mean field spin glass. It is a 
caricature of the Sherrington \&
Kirkpatrick (SK) spin glass model [\rcite{MPV}]. 
Both are spin systems with Ising spins taking value $\pm 1 $. 
In the SK model one has Gaussian pair interactions,  
while in the REM one has Gaussian multibody interactions
of any order.
The Hamiltonian  of the REM is
$$
H(\s)\equiv -\frac{ \sqrt{N}}{2^{N/2}}
\sum_{\a \subset \{1,\dots,N\}}J_{\a}\s_{\a}
\Eq(1.00001)
$$

where the sum is over all the $2^N$ subsets of $\{1,\dots,N\}$, 
$(J_\a, \a \subset \{1,\dots,N\})$
is a family of i.i.d.  standard Gaussian variables defined on a common 
probability space $(\O,\S,\QQ)$ and
$\s_a\equiv \Pi_{i\in \a} \s_i$ with the convention that 
$\s_{\emptyset}=1$. It turns out that the random variables $H(\s)$ and 
$H(\s')$ corresponding to different configurations $\s\not=\s'$ are independent
Gaussian 
variables with zero mean and variance $N$.
The equilibrium statistical mechanics of the REM  has been well  studied,
e.g., in a non rigorous way, in [\rcite{B1},\rcite{B2}] and, in a rigorous
way, in [\rcite{Ei},\rcite{OP},\rcite{GMP}].
We quote some of the (rigorous) results that will be important for understanding the
dynamics.
Given $\b\geq 0$, the inverse temperature,
let us denote by
$$
Z_N\equiv Z_N(\b)= \sum_{\s} e^{-\b H(\s)}
\Eq(1.00002)
$$
the finite volume partition function and by
$$
F_N(\b)= \frac 1{N} \log Z_N(\b)
\Eq(1.00003)
$$
the finite volume free energy.

It was proved in [\rcite{OP}] that for all $\b\geq 0$
the limit $\lim_{N\rightarrow \infty} F_N(\b)=F(\b)$
exists $\QQ$-almost surely  and in $L^p(\O,\S,\QQ)$ for $1\le p<\infty$.  
$F(\b)$ equals $\b^2/2 +\b_c^2/2$ for $\b<\b_c$ and $\b_c\b$ for
$\b\ge \b_c$, as expected from the results of [\rcite{B1}]. 
$F(\b)$ is therefore a non random function which is twice differentiable in $\b$
but the second derivative has a jump at  $\b_c=\sqrt {2\log 2}$. 
This is called in the physics literature a third order phase transition.
Another important fact is that,   
depending wether we are in a high temperature regime ($\b<\b_c$) 
or in a low temperature one ($\b\geq \b_c$), 
not only does the free energy change from a quadratic function of $\b$
to a linear one but the difference between the finite volume free energy  
and its infinite volume limit 
is exponentially small in $N$ in the high temperature case, 
whereas, in the low temperature regime, $F_N(\b)-F(\b)$ behaves as
$C(\omega,\b,N)\frac {\log N}{N}$, for some random function $C(\omega,\b,N)$.
$C(\omega,\b,N)$ converges in $\QQ$-probability to a non-random limit but 
does not converge $\QQ$-almost surely and  the $\QQ$ almost-sure cluster 
set of  $C(\omega,\b,N)$ was identified in [\rcite{GMP}].

Let us now discuss the dynamical properties of the model. We consider the Metropolis 
dynamics. (See \eqv(4.4)).  
A  first step in the study of the dynamics for the REM was done 
in [\rcite{FIKP}]. 
There the spectral gap, $\lambda_N$ of the usual single spin flip
metropolis dynamics in volume
$N$ is studied. In particular it was proved that  for {\it all} 
inverse temperatures $\b>0$
we have, $\QQ$-almost surely 
$$
\lim_{N\uparrow\infty} - \frac 1N \log \l_N = \b \b_c
\Eq(1.0001)
$$

Moreover $\QQ$-almost sure finite size corrections are also given 
in [\rcite{FIKP}]: we have
$$
\b\b_c-c\b \sqrt{\frac{\log N}{N}} \le - \frac 1N \log \l_N
\le \b\b_c + c \b  \sqrt{\frac{\log N}{N}}
\Eq(1.0002)
$$
$\QQ$-almost-surely, for all but a finite number of indices $N$, 
for some constant $c$.

However  one would have expected the dynamics to 
present a kind of  transition  as the previously mentioned {\it static} phase 
transition that can be  seen on the free energy $F(\b)$. 
Such a {\it dynamical} transition is 
not seen on the spectral gap.

Thus we are lead to the following question: how can we see a dynamical phase transition 
on the single spin flip dynamics ? 

The inverse spectral gap can be used as an estimate for the thermalization time 
of the dynamics. For the Metropolis dynamics,  $1/\l_N$ is actually a sharp upper bound for 
the time it takes for the
dynamics to reach equilibrium, whatever was the initial law. In particular we may
consider the dynamics issued from a given configuration. 
The REM is rather pathological in the sense that  
the configurations of lowest energy ( of order $-\b_c N$) are surrounded 
(in a sense of a single spin flip)
by configurations of energy of order at most $ \pm \sqrt{N \log N}$.  
The   bounds in \eqv(1.0002) follow from this fact. Starting the dynamics at a 
configuration of lowest energy, we have to wait for a time of order  $e^{N\b\b_c}$ 
before the first spin flip. 
As we see, the time to reach equilibrium starting from a configuration of minimal
energy is therefore of order $e^{N\b\b_c}$. 

In the low temperature regime, $\b>\b_c$, the equilibrium measure is concentrated on
these configurations of minimal energy. But in the high temperature regime, $\b<\b_c$, 
the invariant measure does not charge too much  
these configurations with minimal energy. In fact 
the invariant measure
has its mass concentrated on configurations with energy of order 
$-\b N$. This follows 
from results in  [\rcite{OP}]. In a certain sense, when $\b<\b_c$, 
it is therefore 'un-natural' 
to compute the thermalization time starting from a configuration of minimal energy.

We shall therefore change our point of view: instead of considering any initial law, 
we shall rather estimate the time to equilibrium when the dynamics 
starts from the uniform probability. Doing this we expect the dynamics to avoid 
the configurations of minimal energy ( in the high temperature regime), and thus we 
hope to see a dynamical phase transition. 

Using generalized Poincar\'e inequalities, we get upper bounds for the time to 
equilibrium starting from the uniform law, say $T_N$. We prove that, when $\b<\b_c$, 
then  
$$\limsup \frac 1 N \log T_N\le \b^2\Eq(1.22)$$
Comparing \eqv(1.22) with \eqv(1.0002), one sees that the thermalisation time is 
much shorter than the inverse spectral gap. In other words, in the high temperature
regime, starting from the uniform law, the dynamics reaches equilibrium much faster
than starting from one of the configurations of minimal energy. 
These results can be interpreted as a first step towards a proof of the existence of a
dynamical phase transition. Actually we expect \eqv(1.22) to be sharp i.e. we expect
$\frac 1N \log T_N$ to converge to $\b^2$, for all $\b<\b_c$. In the low temperature
regime, the asymptotics of $T_N$ should be given by the inverse spectral gap i.e. one
expects $\frac 1N \log T_N$ to converge to $-\b\b_c$ for all $\b\ge\b_c$. Thus one
would see the dynamical phase transition for the Metropolis dynamics.

Remember that the  Hamiltonians $H(\s), \s \in \{-1,+1\}^N$ of the REM 
form a family of i.i.d, Gaussian random variables 
with mean zero and variance $N$, 
defined on 
some probability space, say  $({\O},\Sigma,\QQ)$.
Given $\b\geq 0$, the inverse temperature, 
the Gibbs measure is defined by 
$$
\p_{\b}(\s)\equiv \frac{e^{-\b H(\s)}}{Z_N}
\Eq(4.2)
$$
where $Z_N$ is defined in \eqv(1.00002).
For a given realization of the Hamiltonian,
we consider the Metropolis dynamics, $X(t)=X_N(t)$:  
 $X(t)$ is the continuous time Markov process 
defined on ${\XX} \equiv \{-1,+1\}^N$ by the transition rates: 
$$
P(\s,\s')=
\cases{\frac 1N\exp\{-\b(H(\s')-H(\s))^+\} &{if} $ ||\s'-\s||=1$\cr
0 & {if} $ ||\s'-\s||>1$\cr
}\Eq(4.4)
$$
where $a^+=\max\{a,0\}$ and $||x||=\frac{1}{2}\sum_{i=1}^{N}|x_i|$.
 $\p_{\b}$ is invariant, ergodic, and reversible for this dynamics. 

 The associated Dirichlet form on $L_2({\XX},\pi_\b)$ is given by
$$
\EE(f,g)=\frac 1{2N Z_N(\b)}
\sum_{x,y} (f(x)-f(y))(g(x)-g(y))e^{-\b (H(x)\vee H(y))}
\Eq(4.5)
$$
With the notation of part II, 
$$
Q(e)=\frac 1{NZ_N(\b)} e^{-\b (H(x)\vee H(y))}
$$ for $e=(x,y)$ with 
$\Vert x-y\Vert=1$. 

From \eqv(1.0002), one deduces that for any fixed initial law $\h$, 
and any $\g>\b\b_c$, 
then, $\QQ$.a.s. 
$$d_{TV}(\LL_\h(X(e^{\g N})),\p_\b)\rightarrow 0\Eq(4.6.1)$$ 

From now on, we assume that $\b\le\b_c$. 
Given a probability measure $\eta$ on $\XX$, and $t\in \R$, let
$\LL_{\eta}(X(t))$
be the law of the process at time $t$ starting from the initial measure $\eta$.

Given $\e >0$,  $c>0$ and a probability measure $\eta$,  we define the time 
$ T_N(\e,c,\eta)$ to
reach equilibrium
starting from $\eta$, up to
$\e$ on a subset of $\QQ$-probability greater than $1-e^{-c N}$ by 
$$
T_N(\e,c,\eta)\equiv \inf \left\{T\ge 0 : \QQ\left[\sup_{s\ge T} 
d_{TV}\left(L_{\eta}(X(s)),\p_\b\right) \le \e \right]\ge 1-e^{-c N}\right\}
\Eq(4.8)
$$

The main result of this section is

{\bf \Theorem(theo5)} {\it Let $\eta$ be the uniform probability 
measure on $\XX$. 
Then  for all $c>0$, $\e>0$ and for all $\b \le \b_c$
$$
\limsup_{N \uparrow \infty} \frac 1N \log T_N(\e,c,\eta)
\le \b^2
\Eq(4.9)
$$
}
We can also prove estimates when $\eps$ goes 
to $0$ as $N\rightarrow \infty$. We consider two cases: 
$\eps$ going to $0$ polynomialy and as a  stretched exponential. 

{\bf \Theorem(theo6)} {\it Let $\eta$ be the 
uniform probability measure on $\XX$. 
There exists a constant $c_1>0$, such that 
for all $c>0$, there exists a constant 
$C_0=C_0(c,\b)$ such that 
$$\eqalign{
&\frac 1N \log T_N(e^{- N^{1/4}(\log N)^{3/4}}, c,\eta)\cr
\le &\b^2 +  2\b \b_c \left({ \frac{c_1(1+c)\log N}{N}}\right)^{1/2} + c_2(\b,c)
\left( \frac {\log N}{N}\right)^{1/4}+C_0 (\frac{\log N}N)^{3/4}
}\Eq(4.901)
$$
where
$$
c_2(\b,c)\equiv \b\left(12 \b\b_c \sqrt{ c_1(1+c)}\right)^{1/2} 
+ \frac 14\frac{ \b^2+\b_c^2}{\b \b_c \sqrt{ c_1(1+c)}}
\Eq(4.9010)
$$
Moreover  for all $\d>0$
$$
\eqalign{
&\frac 1N \log T_N(N^{-\d}, c,\eta)
\le \b^2 + 
\b\left( 12 \b\b_c\sqrt{c_1(1+c)}\right)^{1/2} (\frac{\log N}N)^{1/4}\cr
& +2\b \b_c\left({ \frac{c_1(1+c)\log N}{N}}\right)^{1/2} 
+\frac 14\frac{\b^2+\b_c^2}{\b\b_c}\frac\d{\sqrt{c_1(1+c)}}(\frac{\log N}N)^{1/2}
+C_0 \frac{\log N}N
}\Eq(4.902)
$$
}

As a corollary, we get

{\bf \Corollary(cor7)}{ \it Let $\eta$ be the 
uniform probability measure on $\XX$. 
For all $\g >\b^2$,
with a $\QQ$-probability 1, for all but a finite number of indices $N$
$$  
d_{TV}\left( \LL_{\eta}(X(e^{\g N})), \p_{\b}) \right) 
\le e^{-N^{1/4} (\log N)^{3/4}}
\Eq(4.11)
$$
Moreover for all $\d>0$, if 
$$
t_N= \exp[ \b^2 N+ \sqrt{12\b\b_c} N^{3/4}(\d \log N)^{1/4}+ 
(2\b\b_c + \d (\b^2+\b_c^2)) (c_1 N\log N)^{1/2}]
\Eq(4.110)
$$
then
$$  
d_{TV}\left( \LL_{\eta}(X(t_N)), \p_{\b}) \right) \le \frac{1}{N^{\d}}
\Eq(4.111)
$$
}

 The error terms in the bound \eqv(4.901) have no reason to be
optimal. However in \eqv(1.0002) the order of magnitude of the 
error terms are optimal as it was observed in [\rcite{FIKP}]. 
\smallskip

To prove the theorems we will need estimates for 
the constants, $\LL_{\pi_\b}(p)$ and
$\LL_\eta(p)$ using \eqv(2.20). This will be done now and the result
will be collected in the Proposition \eqv(prop13).

\smallskip

These estimates will also depend on the choice of paths $\g(x,y)$. To 
estimate the spectral gap, the following set of paths was introduced in  
[\rcite{FIKP}] and they work also here:
given $i\in\{1,\dots,N\}$, and $x,y\in \XX$, such that $x_i\neq y_i$  
let $\g^i(x,y)$ be the path starting at $x$ and ending at $y$
obtained by flipping  the disagreeing spins, starting at the site $i$ and then
going cyclically. Let $\G^i=\{\g^i(x,y), x,y\in \XX\}$. Given 
$x,y$ and $\g(x,y)$, let $\overline{\g(x,y)}$ be the set of points visited by the
path and $\g^o(x,y)=\overline{\g(x,y)}\setminus \{x,y\}$ the set of the
interior points of the path. Note  that if the number of
discrepancies between $x$ and $y$ is $n$ then there exist $n$ interior
disjoint paths in $\{\g^i(x,y), i=1,\dots,N\}$. This comes from the fact that if
$i_1,\dots,i_n$ are the $n$ sites where $x$ and $y$ disagree, then the paths
$\g^{i_1}(x,y),\dots,\g^{i_n}(x,y)$ are interior disjoint. 
The set of paths we will construct will depend on the realization of $H(x)$: 
it is a random set. 
Given a positive number $c_e$, we will say that a point $z$ is good 
if $H(z) \le \sqrt{(1+c_e)2N\log N}$. Call $G$ the set of good points. 
If $z$ is not good, we call it 'bad' and write  $B$ for the set of bad points. 
A path is good if all its interior points are good. 
Note that  we need to select a path for any pair 
of points $(x,y)$, and the typical number of bad points is of order
$2^N2^{-(1+c_e)\log N}$. We cannot neglect  good paths $\g(x,y)$
with  bad end points $x$ or  $y$ or both.
We construct the set of paths $\G$ according to the following rules:

For $||x-y||_1\ge \frac {N}{\log N}$, if there is a good path in  
$\{\g^i(x,y), i=1,\dots,N\}$, choose the first and put
it in $\G$  ;
otherwise, choose $\g^1(x,y)$.
For $||x-y||_1<\frac {N}{\log N}$, if there exists a good site 
$z$ in $\XX$ such that 
$||x-z||\geq\frac {N}{\log N}$,
$||y-z||\geq\frac {N}{\log N}$
and if there are good paths, one in  
$\{\g^i(x,z), i=1,\dots,N\}$
 and another in
$\{\g^i(z,y), i=1,\dots,N\}$
such that the union of
these two good paths is a self avoiding path, then we select this union 
as the path connecting $x$ and $y$
in $\G$. (Note that this is a good path since $z$ is good); otherwise, select
$\g^1(x,y)$. Note that all the paths constructed in this way  have length smaller
than $N$.
A fundamental result that can be easily proven by keeping the $\QQ$-probability
in the proof of  the proposition 4 .1 in [\rcite{FIKP}] is

{\bf \Proposition(prop8)}{\it For all $c_e>0$, 
there exists $N_0(c_e)$ such that for all $N\ge N_0(c_e)$, 
with a $\QQ$-probability $\geq 1-e^{-c_e N}$, 
all the paths of the previous set $\G$ are good i.e. 
satisfy $H(z)\le \sqrt{(1+c_e) N\log N}$, for all $z\in
\g(x,y)\setminus (x,y)$, for all $(x,y)\in \XX^2$. Moreover they have
a  length smaller than N.}  
\smallskip
We say that an edge
$e=(x,x')$ is good,  if $x$ and $y$ are good, this will be denoted by $e\in
\GG$, otherwise the edge is bad: $e\in \BB$.
Note the important fact that, with our construction, a
 given  edge  $e=(x,x')$ belonging to $\G$ 
can have at most one bad point among $x$ and  $x'$. 

\medskip

Let us first estimate, $\LL_{\pi_\b}(p)$, see \eqv(2.20).
The weights $\l(x)$ are chosen in the following way:
Let $d$ be  such that $\b< d\le \b_c$, to  be chosen later.
We set $d=\b(1+\z)$ with $ 0<\z< (\b_c -\b)/\b$.
Let
$$
\l(x)=\cases{ 1 &if $H(x)\ge -d N$\cr
\l& otherwise.\cr}
\Eq(4.12)
$$
where 
$$
\eqalign{
\l&\equiv \left( \sum_{x\in \XX} e^{-\b H(x)}\1_{\{H(x)\le -d N\}}\right)^\r
\cr
& \equiv \left( Z_N(\b,\le -d)\right)^\r}
\Eq(4.13)
$$
for some $\r >0$ to be chosen later.

First we  consider the  first term in the right hand side of \eqv(2.20)
the other ones will be treated later.  Let us denote

$$
R(d,\r,p)
\equiv \sum_{x\in \XX}\pi_\b(x)(\l(x))^{p/1-p}
\Eq(4.14)
$$

{\bf \Lemma(lem9)}{\it Let $ \z >0$, $0<\r<1$ and $0<p<1/2$ that  satisfy 
$0<\z\le (\b_c-\b)/\b$ and 
$$
\frac {p\r}{\z^2(1-p)} \le \frac 12 \frac {\b^2}{\b^2+\b_c^2}
\Eq(4.15)
$$
There exists an absolute constant 
$c_1$, and, for any $c>0$, 
there exists $N_0(\b,c,\z)$ such that for any $N\ge N_0(\b,c,\z)$ such that
$$
\sqrt{ \frac{N}{(1+c)\log N}} \ge \frac {12 \b_c}{\z^2 \b} \sqrt {c_1}
\Eq(4.151)
$$
then, with a $\QQ$-probability$\ge 1-e^{-cN}$, we have 
$$
R(d,\r,p)
\le 2
\Eq(4.16)
$$
where $d=\b(1+\z)$} 

{\proof}
 Let us denote 
$$
R^1(d,\r,p)
\equiv\sum_{x\in\XX}\pi_\b(x)\1_{\{H(x)\geq -d N\}}
\Eq(4.17)
$$ 
and $R(d,\r,p)\equiv R^1(d,\r,p)+R^2(d,\r,p)$. 
We have $R^1(d,\r,p)\le 1$ and, to estimate $R^2(d,\r,p)$, we use the
following lemma that will be proved in the section  VI.

{\bf \Lemma(lem10)}{\it  There exits a constant $c_1$, such that for 
all $c>0$, there exists a $N_0(\b,c,\z)$ such that for all $N\ge N_0(\b,c,\z)$
with a $\QQ$-probability $\ge 1-e^{-cN}$
$$
Z_N(\b,\le-d)\le 
e^{\b\b_c \sqrt{c_1(1+c) N\log N }}
e^{N[\b d-\sfrac{d^2}{2}+\sfrac{\b_c^2}{2}]}
\Eq(4.18)
$$
 and
$$
Z_N(\b)\ge 
e^{-\b\b_c \sqrt{c_1(1+c) N\log N }}
e^{N[\sfrac{\b^2}{2}+\sfrac{\b_c^2}{2}]}
\Eq(4.19)
$$
}

Note that $$R^2(d,\r,p)=\frac{Z_N(\b,\le-d)^{1+\r p/(1-p)}}{Z_N(\b)}$$
Therefore Lemma \eqv(lem10) implies that 
$$
R^2(d,\r,p)\le 
e^{[2+\sfrac{p\r}{1-p}]\b\b_c \sqrt{c_1(1+c) N\log N }}
e^{-N[\sfrac{(d-\b)^2}{2}]}
e^{N\sfrac{p\r}{1-p}[\b d-\sfrac{d^2}{2}+\sfrac{\b_c^2}{2}]}
\Eq(4.20)
$$
 Now using \eqv(4.15), we get
$$
\frac {\r p}{1-p}\left[ \b d-\frac{d^2}{2} +\frac{\b_c^2}{2}\right]
\le \frac 12 \frac {(d-\b)^2}{2}
\Eq(4.200)
$$
Using 
 $0<\r<1$ and $ 0<p<1/2$, we have $\r p/(1-p) \le 1$
therefore \eqv(4.151) implies 
that 
$$
\left[ 2+ \frac {p\r}{1-p}\right]\b\b_c \sqrt{c_1(1+c)N\log N}
\le \frac 12 \frac {(d-\b)^2}{2} N
\Eq(4.2001)
$$
from which we  immediately get \eqv(4.16).\qed

Now we  estimate the other term in the right hand side of \eqv(2.20). Let us 
denote 
$$
\frac{1}{2\LL^*_{\p_\b}(1)} \equiv
\sum_e \frac 1{Q(e)}\left[\sum_{x,y:\g(x,y)\ni e} \frac {\pi_\b(x)}{\l(x)}
\frac{\pi_\b(y)}{\l(y)}\right]^2
\Eq(4.21)
$$

{\bf \Proposition(prop11)}{\it We assume that  $2(1-\r)<1$. 

There exists a constant $c_1$, such that for all
$c >0$ and all $\z>0$, $\z<(\b_c-\b)/\b$,  satisfying \eqv(4.15),  
there exist $N_0(\b,c,\z)$ such that for all $N\ge N_0(\b,c)$,  
with a $\QQ$-probability $\ge 1-e^{-cN}$, we have 
$$
\frac{1}{2\LL^*_{\pi_\b}(1)} \le  
22 N^4 e^{2\b\b_c \sqrt{c_1(1+c) N\log N}}e^{\b dN}
\Eq(4.22)
$$
}

{\proof}
 Let us write 
$$
\frac{1}{2\LL^*_{\pi_\b}(1)} = L_{\BB} +L_{\GG}
\Eq(4.23)
$$  
where $L_{\GG}$ is the same as \eqv(4.21) but with the sum $\sum_{e}$
restricted to good edges and $L_{\BB}$ with bad edges.

Let us first consider $L_{\BB}$.
Using convexity and symmetry,we can write 
$$
L_{\BB}\le 3L_{\BB}(\ge,\ge)+ 6 L_{\BB}(\ge,<)  \Eq(4.24)
$$ 
where $L_{\BB}(\ge,\ge)$ is the same as in \eqv(4.21) but with the following
restrictions: 
$e\in \BB, H(x)\ge -d N, H(y)\ge -d N$. $L_{\BB}(\ge,<)$ is defined similarly. 

Let $ U\equiv  \{x\ ;\ H(x)\ge -d N\}$ and $ D \equiv \{x\ ;\ H(x)< -d N\}$. 
Since a bad edge is the first or the last edge of the path, 
if $e=(z,z')\in\g(x,y)$ is a bad edge then  
we have either $z\in B$ and $x=z$ or $z'\in B$ and $z'=y$. 
By symmetry it is sufficient 
to consider the first case. Then $1/Q(e)=N Z_N(\b) \exp(\b H(z))$. 
Note in particular that it is 
not possible to have $e=(z,z')\in\g(x,y)$, $e$ bad and both $x$ and $y$ 
belonging to $D$. This is the reason why we do not have a term $L_{\BB}(<,<)$.

$$
\eqalign{
L_{\BB}(\ge,\ge)
&\le \frac{2N}{Z_N^3(\b)} \sum_{e=(z,z')}e^{\b H(z)}
\left[ \sum_{x,y\in U\times U}e^{-\b H(x)} e^{-\b H(y)}\1_{\{x=z\}}\right]^2\cr
&=\frac{2N}{Z^3_N(\b)}\sum_{e=(z,z')}e^{-\b H(z)}
\left[ \sum_{y\in U}e^{-\b H(y)} \right]^2 \le 2N\cr
}\Eq(4.25)
$$

Using similar arguments and recalling \eqv(4.13), we get

$$
\eqalign{
L_{\BB}(\ge,<)&\le \frac{2N}{Z_N^3(\b)} \sum_{e=(z,z')}e^{\b H(z)}
\left[ \sum_{x,y\in U\times D}
e^{-\b H(x)} \frac{e^{-\b H(y)}}{Z^{\r}_N(\b,\le -d)}\1_{\{x=z\}}\right]^2\cr
&\le \frac{2N}{Z_N^3(\b)}\sum_{e=(z,z')}e^{-\b H(z)}
\left[ \sum_{y\in D} \frac{e^{-\b H(y)}}{Z^{\r}_N(\b,\le -d)}\right]^2\cr
&\le \frac {2N}{Z_N^2(\b)} \left[Z_N(\b,\le -d)\right]^{2(1-\r)}\cr
}
\Eq(4.26)
$$

Using \eqv(4.18),  \eqv(4.19), and $2(1-\r)\le 1$ we get
$$
\eqalign{
L_{\BB}(\ge,<) &\le 2N 
e^{3\b\b_c \sqrt{c_1(1+c) N\log N}}
e^{-N[\sfrac{\b^2+\b_c^2+(d-\b)^2}{2}]}\cr
&\le 2N e^{-N[\sfrac{\b^2+\b_c^2}{2}]}\cr
&\le 2N\cr}
\Eq(4.27)
$$
where  we have used \eqv(4.151) at the second step.
We have proved that 
$$
L_{\BB}\le 18 N
\Eq(4.28)
$$

We consider now $L_{\GG}$ As before, using convexity and symmetry, we write
$$
L_{\GG}\le 4  L_{\GG}(\ge,\ge) + 8 L_{\GG}(\ge,<) +4L_{\GG}(<,<)
\Eq(4.29)
$$ 

We first consider $ L_{\GG}(<,<)$. Since for a good edge $e=(z,z')$, 
we have $H(z)\vee H(z')\le 
\sqrt{(1+c_e)N\log N}$, we therefore get

$$
\eqalign{
L_{\GG}(<,<)&\le \frac{N e^{\b\sqrt{(1+c_e)N\log N}}}{Z^3_N(\b)}
\sum_{e\in \GG}
\left[ \sum_{x,y\in D\times D}
\1_{\{\g(x,y)\ni e\}}\frac{e^{-\b H(x)}}{Z^{\r}_N(\b,\le -d)}
\frac{e^{-\b H(y)}}{Z^{\r}_N(\b,\le -d)}
\right]^2\cr
}
\Eq(4.30)
$$

On  one hand we have
$$
\left[ \sum_{x,y\in D\times D}\1_{\{\g(x,y)\ni e\}}
\frac{e^{-\b H(x)}}{Z^{\r}_N(\b,\le -d)}
\frac{e^{-\b H(y)}}{Z^{\r}_N(\b,\le -d)}
\right]
\le Z_N^{2(1-\r)}(\b,\le-d)
\Eq(4.31)
$$
On the other hand we have
$$
\eqalign{
& \sum_{e\in \GG}
 \sum_{x,y\in D\times D}
\1_{\{\g(x,y)\ni e\}}\frac{e^{-\b H(x)}}{Z^{\r}_N(\b,\le -d)}
\frac{e^{-\b H(y)}}{Z^{\r}_N(\b,\le -d)}\cr
& = 
 \sum_{x,y\in D\times D}
\frac{e^{-\b H(x)}}{Z^{\r}_N(\b,\le -d)}
\frac{e^{-\b H(y)}}{Z^{\r}_N(\b,\le -d)}
\sum_{e\in \GG}
\1_{\{\g(x,y)\ni e\}}\cr
&\le N Z_N^{2(1-\r)}(\b,\le-d)\cr
}\Eq(4.32)
$$
where at the last step we have used that the length of a path is smaller than
$N$.
Therefore using \eqv(4.31) and \eqv(4.32) in \eqv(4.30), then \eqv(4.18) and
\eqv(4.19) and at last $4(1-\r)\le 2\le 3$ 
we get
$$
\eqalign{
L_{\GG}(<,<) &\le N^2  e^{\b\sqrt{(1+c_e)N\log N}} 
\frac { Z_N^{4(1-\r)}(\b,\le-d)}{ Z^3_N(\b)}\cr
&\le N^2 e^{\b\sqrt{(1+c_e)N\log N}}
e^{6\b\b_c \sqrt{c_1(1+c) N\log N}}
e^{-N\sfrac{3(d-\b)^2}{2}}\cr
&\le N^2 e^{7\b\b_c \sqrt{ c_1(1+c)N \log N}}
e^{-N\sfrac{3(d-\b)^2}{2}}\le N^2\cr
}\Eq(4.33)
$$
Where we have used $\b_c=\sqrt{2\log 2}>1$, \eqv(4.151) and we
have chosen $c_e=c$ and $c_1 >1$.

Consider now $L_{\GG}(\ge,<)$. Using exactly the same kind of arguments,
using \eqv(4.18) and \eqv(4.19), and $2(1-\r)\le 1$
we get
$$
\eqalign{
 L_{\GG}(\ge,<)&\le N^2 e^{\b\sqrt{c_1(1+c)N\log N}} 
\frac { Z_N^{2(1-\r)}(\d,\le-d)}{ Z_N(\b)}\cr
&\le 
 N^2 e^{3\b\b_c \sqrt{c_1(1+c)N\log N}} 
e^{-N\sfrac{(d-\b)^2}{2}}\le N^2\cr
}\Eq(4.34)
$$
where at the last step we have used \eqv(4.151).

We consider now $L_{\GG}(\ge,\ge)$. Since the edge is good, we have
$$
\eqalign{
L_{\GG}(\ge,\ge)&\le \frac{Ne^{\b\sqrt{(1+c_e)N\log N}}}{Z^3_N(\b)}
\sum_{e\in \GG}
\left[ \sum_{x,y\in U\times U}
\1_{\{\g(x,y)\ni e\}}e^{-\b H(x)} e^{-\b H(y)}\right]^2\cr
&\le {Ne^{\b\sqrt{(1+c_e)N\log N}}} \sup_{e}
\left[ \sum_{x,y\in U\times U}\1_{\{\g(x,y)\ni e\}}
\frac{e^{-\b H(x)} e^{-\b H(y)}}{Z_N(\b)}\right]\cr
&\quad\quad \times 
\sum_{e\in \GG}
\left[ \sum_{x,y\in U\times U}\1_{\{\g(x,y)\ni e\}}
\frac{e^{-\b H(x)} e^{-\b H(y)}}{Z^2_N(\b)}\right]\cr
&\le {N^2e^{\b\sqrt{(1+c_e)N\log N}}} \sup_{e}
\left[ \sum_{x,y\in U\times U}\1_{\{\g(x,y)\ni e\}}
\frac{e^{-\b H(x)} e^{-\b H(y)}}{Z_N(\b)}\right]\cr
}\Eq(4.35)
$$

To continue we will need an adaptation of [\rcite{FIKP}].
Let us call
$$
\L{(d)}\equiv \sup_{e}
\left[ \sum_{x,y\in U\times U}\1_{\{\g(x,y)\ni e\}}
\frac{e^{-\b H(x)} e^{-\b H(y)}}{Z_N(\b)}\right]
\Eq(4.36)
$$
 
Recalling that the paths in $\G$ are constructed using paths in
$\cup_{i=1}^N\G^i$, we get immediately
$$
\L{(d)}\le N \sup_{1\le i \le N }
\L^{(i)}({d})
\Eq(4.37)
$$
where  $\L^{(i)}(d)$ is as in \eqv(4.36) but  with  paths in $\G^i$.
It is enough to consider the case $i=1$ the other ones being similar.
Now for a given edge $e=(z,z')$, there exists a $j\in \{1,\dots,N\}$
such that $z'=z^j$, that is $z'$ is the configuration obtained from $z$
by flipping the spin at the site $j$. Note at this point that the set
of all $(x,y): \g(x,y)\ni e$ for $\g \in \G^1$ is exactly
$$
\bigcup_{x\in \{-1,+1\}^{j-1}}
\bigcup_{ y\in \{-1,+1\}^{N-j}}
((x_1,\dots,x_{j-1},z_{j},\dots,z_{N}),(z_1,\dots,z_{j-1},-z_{j}
,y_{j+1},\dots,y_{N})) 
\Eq(4.38)
$$
Denoting $z_{>j}\equiv(z_{j+1},\dots,z_N)$,  
$z_{<j}\equiv(z_1,\dots,z_{j-1})$,
$$
Z^{(1)}_{j-1}(\b,\ge -d)[z_{j},z_{>j}]
\equiv
\sum_{x\in \{-1,+1\}^{j-1}} e^{-\b H(x,z_j,z_{>j})}
\1_{\{H(x,z_j,z_{>j})\geq -d N\}}
\Eq(4.39)
$$
and
$$
Z^{(1)}_{N-j}(\b,\ge -d)[z_{<j},-z_{j}]
\equiv
\sum_{y\in \{-1,+1\}^{N-j}} e^{-\b H(z_{<j},-z_j,y)}
\1_{\{H(z_{<j},-z_j,y)\geq -d N\}}
\Eq(4.40)
$$
we get immediately:
$$
\L^{(1)}(d)\le
 \sup_{z\in \XX}
\frac {1}{Z_N(\b)}
Z^{(1)}_{j-1}(\b,\ge -d)[z_{j},z_{>j}]
Z^{(1)}_{N-j}(\b,\ge -d)[z_{<j},-z_{j}]
\Eq(4.41)
$$

To continue, we need the following lemma that will be proved in the
 next section. 

{\bf \Lemma(lem12)}{\it  There exists a constant $c_1>0$, such that for all
$c>0$, if $c_u=c_1^{-1}( 2\log 2+c)$, 
then we can find an $N_0=N_0(\b,c)$ such that
for all $N>N_0(\b,c)$, with a $\QQ$-probability
$\geq 1-e^{-cN}$, if  we call $M\equiv\sqrt{ N/\log_2(Nc_u)}$ and 
$j-1\equiv \a N$ then
$$
\sup_{z_j,z_{>j}}Z^{(1)}_{j-1}(\b,\ge -d)[z_{j},z_{>j}]
\le \sqrt{N}e^{\b\b_c \sqrt{N\log_2(c_u N)}}
\left( 2^j + \ZZ_N(\b,d,\a)\right)
\Eq(4.42)
$$
where}
$$
\ZZ_N(\b,d,\a)=
\cases{ e^{\b dN} & if $\sqrt{\a M^2-1} <\frac \b{\b_c} M$\cr
 e^{\b dN} + e^{N[\sfrac{\b^2}{2}+\a \sfrac{\b_c^2}{2}]}&if
$\frac \b{\b_c} M\le \sqrt{\a M^2 -1} < \frac{d}{\b_c} M $\cr
e^{N[\sfrac{\b^2}{2}+\a \sfrac{\b_c^2}{2}]}&if
$\frac d{\b_c} M\le \sqrt{\a M^2 -1}$\cr
}\Eq(4.43)
$$

Now inserting \eqv(4.43) for $j-1=\a N$ and $N-j=(1-\a)N$ in 
\eqv(4.41), considering the nine resulting terms, using
$$
\b d < \b\b_c <\frac {\b^2}{2} +\frac {\b_c^2}{2}
\Eq(4.44)
$$
to simplify the computations,
and maximizing over $\a\in [0,1]$, it is just a long task to get
$$
\eqalign{
\L^{(1)}(d)&\le
\sqrt {N}e^{\b\b_c \sqrt{N\log_2(c_u N)}}e^{\b dN}\cr
&\le \sqrt{N} e^{\b\b_c \sqrt{c_1(1+c)N\log  N}}e^{\b dN}\cr
}\Eq(4.45)
$$
with a $\QQ$-probability $\geq 1-e^{-cN}$.  
Inserting  \eqv(4.45) and \eqv(4.37) in \eqv(4.35),  we get
$$
L_{\GG}(\geq,\geq) \le
{N^4 e^{2\b\b_c\sqrt{c_1(1+c)N\log N}}}  
e^{\b d N}
\Eq(4.46)
$$
 Using \eqv(4.23),\eqv(4.28),\eqv(4.33),\eqv(4.34) and \eqv(4.46)
we get \eqv(4.22).\qed

Now we estimate $\LL_{\eta}(p)$ see \eqv(2.20),  when $\eta$ is the uniform
measure on $\XX$. We take the weights $\mu(x) =1$.
Since we have already estimated the first factor in Lemma \eqv(lem9),
 it remains to estimate

$$
\frac {1}{\LL^*_{\eta}(1)}\equiv
\sum_{e} \frac 1{Q(e)} 
\left[ \sum_{x,y:\g(x,y) \ni e} 
\frac{\pi_{\b}(x)}{\l(x)} \frac 1{2^N}\right]^2
\Eq(4.47)
$$

{\bf \Proposition(prop13)}{\it 
There exists a constant $c_1$, such that for all
$c >0$,   for all $\z>0$, $\z <(\b_c-\b)/\b$, there exists $N_0(\b,c,\z)$, such that for all 
$N$  that satisfy \eqv(4.151) and are larger than $N_0(\b,c,\z)$, 
with a $\QQ$-probability $\ge 1-e^{-cN}$, 
$$
\frac {1}{2\LL^*_{\eta}(1)}\le 
4 N^2 e^{4\b\b_c \sqrt{c_1(1+c)N\log N}}
 e^{\b dN }
\Eq(4.471)
$$}

{\proof}
As before by considering separately the cases where $e \in \GG$ and 
$e\in \BB$, we write
$$
\frac {1}{2\LL^*_{\eta}(1)}\equiv L_{\BB}(\eta)+L_{\GG}(\eta)
\Eq(4.48)
$$
Distinguishing bad and good edges and separating the cases $x\in D$ or $x\in U$,
 we get four terms that we 
call, $L_{\BB}(\eta,\geq)$, $L_{\BB}(\eta, <)$,
$L_{\GG}(\eta,\geq)$ and $L_{\GG}(\eta,<)$.

Let us start with $L_{\BB}(\eta,<)$. We should then have $y=z'\in B$. 
Therefore 
$$
\eqalign{
L_{\BB}(\eta,<)&\le \frac {2N}{Z_N(\b)}
\sum_{e=(z,z')}e^{\b H(z')} 
\left[ \sum_{x\in D,y} 
\frac{\pi_{\b}(x)}{Z^\r_N(\b,\le -d)} \frac 1{2^N}\1_{\{y=z'\}}\right]^2\cr
&\le \frac {2N}{Z_N(\b)}
\sum_{z'}e^{\b H(z')} 
\left[\frac {Z^{(1-\r)}_N(\b,\le -d)}{2^N}\right]^2\cr
&=\frac {2N}{Z_N(\b)}Z_N(-\b)
\frac {Z^{2(1-\r)}_N(\b,\le -d)}{2^{2N}}\cr
}\Eq(4.49)
$$
 Now since it is clear that $Z_N(\b)$ and $Z_N(-\b)$ have the same
distribution and therefore satisfy the same estimates, we get that 
with a $\QQ$-probability $\geq 1-e^{-cN}$,  
$$
\frac {Z_N(-\b)}{Z_N(\b)} \le e^{2\b \b_c \sqrt{c_1(1+c) N \log N}}
\Eq(4.50)
$$
Using now  \eqv(4.18), and $2(1-\r)<1$ we get

$$
L_{\BB}(\eta,<)\le 2Ne^{3\b \b_c \sqrt{c_1(1+c) N \log N}}
e^{N[2\b d-d^2]}
\Eq(4.51)
$$

Consider now $L_{\BB}(\eta,\ge)$. Then $e$ is bad and $x\in U$. We have 
to deal separately with, case 1, $x=z\in B$and, case 2, $y=z'\in B$. 
By convexity
$$
L_{\BB}(\eta,\ge)\le 2 
L_{\BB}(\eta,\ge,1)+ 
2L_{\BB}(\eta,\ge,2) 
\Eq(4.52)
$$
On the one hand we  have
$$
\eqalign{
L_{\BB}(\eta,\ge,1)&\le \frac {N}{Z_N(\b)} \sum_{z}
 e^{-\b H(z)}\left[ \sum_{y} \frac {1}{2^N}\right]^2\cr
&=N\cr
}\Eq(4.53)
$$
On the other hand we have
$$
\eqalign{
L_{\BB}(\eta,\ge,2)&\le \frac {N}{Z_N(\b)} 
\sum_{z'}
e^{-\b H(z')}\left[ \sum_{x\in U}e^{-\b H(x)} \frac 1{2^N}
\right]^2\cr
&\le \frac {N}{Z_N(\b)} Z_N(-\b) \left[ \frac
{Z_N(\b)}{2^N}\right]^2\cr
&\le  Ne^{4\b \b_c \sqrt{c_1(1+c) N \log  N}}e^{\b^2 N}\cr
}\Eq(4.54)
$$

 Collecting \eqv(4.51), \eqv(4.53) and \eqv(4.54), we get
$$
\eqalign{
L_{\BB}(\eta) &\le 2N e^{4\b \b_c \sqrt{c_1(1+c) N \log  N}}
\left(e^{\b^2 N} + e^{[2\b d-d^2]N}\right)\cr
&\le 4N e^{4\b \b_c \sqrt{c_1(1+c) N \log N}}
e^{\b d N}\cr
}\Eq(4.55)
$$
where, at the last step, we used that $ d \ge \b$ and therefore 
$2\b d-d^2 \le \b d$.

Consider now $L_{\GG}(\eta,<)$. Since we consider now good edges, we have 
$$
\eqalign{
L_{\GG}(\eta,<)&\le \frac{ Ne^{\b\sqrt{(1+c_e)N\log N}}}{Z_N(\b)}
\sum_{e\in \GG}\left[\sum_{x\in D,y}\frac{e^{-\b H(x)}}
{Z_N^{\r}(\b,\le-d)}\frac 1{2^N}\right]^2\cr
&\le \frac{ Ne^{\b\sqrt{(1+c_e)N\log N}}}{Z_N(\b)}
Z_N^{(1-\r)}(\b\le -d)
\sum_{e\in \GG}\left[\sum_{x\in D,y}\frac{e^{-\b H(x)}}
{Z_N^{\r}(\b,\le-d)}\frac 1{2^N}\right]\cr
&\le\frac{ N^2e^{\b\sqrt{(1+c_e)N\log N}}}{Z_N(\b)}
Z_N^{2(1-\r)}(\b\le -d)\cr
&\le N^2e^{3\b\b_c \sqrt{c_1(1+c)N\log N}}e^{-\frac 12(d-\b)^2N}\cr
&\le N^2}\Eq(4.56)
$$ 
 where we used that  $2(1-\r)<1$ and \eqv(4.151) at the last step.

It remains to consider $L_{\GG}(\eta,\ge)$. Using the fact that 
$e$ is good, we get
$$
\eqalign{
L_{\GG}(\eta,\ge) &\le \frac{ Ne^{\b\sqrt{(1+c_e)N\log N}}}{Z_N(\b)}
\sum_{e\in \GG}\left[\sum_{x\in U,y,\g(x,y)\ni e}e^{-\b H(x)} \frac
1{2^N}\right]^2\cr
&\le Ne^{\b\sqrt{(1+c_e)N\log N}}\sup_{e\in \GG}\left[\sum_
{x\in U,y,\g(x,y)\ni e}e^{-\b H(x)} \frac 1{2^N}\right]\cr
&\quad\quad \times \frac 1{Z_N(\b)}\sum_{e\in \GG}\left[\sum_
{x\in U,y,\g(x,y)\ni e}e^{-\b H(x)} \frac1{2^N}\right]\cr
&\le  N^2e^{\b\sqrt{(1+c_e)N\log N}}\sup_{e\in \GG}
\left[\sum_{x\in U,y,\g(x,y)\ni e}e^{-\b H(x)} \frac 1{2^N}\right]\cr
}\Eq(4.57)
$$

To estimate this last supremum, we use a similar argument as the one
we used to treat \eqv(4.41). Using \eqv(4.42), and the same notation 
as in \eqv(4.39), after a not too long computation, we get that 
with a $\QQ$-probability $\ge 1-e^{-c N}$
$$
\eqalign{
\L(\eta,d)&\equiv
\sup_{e\in \GG}
\left[\sum_{x\in U,y,\g(x,y)\ni e}e^{-\b H(x)} \frac 1{2^N}\right]\cr 
&\le \sup_{1\le j\le N}\sup_{z\in \XX}
Z_{j-1}(\b,\ge -d)[z_j,z_{>j}]2^{-j}\cr
&\le e^{\b\b_c \sqrt{N\log_2(c_u N)}} e^{\b d N}\cr
}\Eq(4.58)
$$
Collecting \eqv(4.57) and \eqv(4.58), we get
$$
L_{\GG}(\eta,\ge) \le N^2e^{2\b\b_c \sqrt{c_1(1+c)N\log N}} e^{\b d N} 
\Eq(4.60)
$$
 Collecting \eqv(4.55), \eqv(4.56) and \eqv(4.60), this entails
\eqv(4.471).\qed 

\smallskip

Now we put together all the results concerning the quantities
$\LL_{\eta}(p)$ and $\LL_{\pi_\b}$. That is collecting Lemmata
\eqv(lem9)  and Proposition\eqv(prop11) and \eqv(prop13), recalling \eqv(2.20)
we have 

{\bf \Proposition(prop14) }
{\it Let $\b<\b_c$, $0<\z <(\b_c-\b)/\b$ and $0<p< 1/2$ satisfy 
$$
\frac {p}{\z^2(1-p)} <  \frac {\b^2}{\b^2+\b_c^2}
$$
There exists an absolute constant $c_1$, such that for all
$c>0$, there exists a
$N_0(\b,c,\z)$ such that for all $N\ge N_0(\b,c,\z)$ 
and $N$ satisfying condition \eqv(4.151)  
then, with a 
$\QQ$--probability $\ge 1-e^{-cN}$, we have 
$$
\frac{1}{\LL_{\pi_\b}(p)}
\le
4^{\sfrac{1-3p +p^2}{p(1-p)}}  
22 N^4 e^{2\b\b_c \sqrt{c_1(1+c) N\log N}}e^{\b^2(1+\z)N}
\Eq(4.571)
$$
and
$$
\frac{ 1}{\LL_{\eta}(p)}
\le
4^{\sfrac{2-3p+2p^2}{p(1-p)}}
4 N^2 e^{4\b\b_c \sqrt{c_1(1+c)N\log N}}
 e^{\b^2(1+\z)N }
\Eq(4.572)
$$
}

{\bf Remark}: the aim of this remark is to discuss the implications of Proposition \eqv(prop14) 
as far as the behaviour of the eigenvectors of the Metropolis dynamics are concerned. To simplify things, 
we only consider the almost sure asymptotics of the first non trivial eigenvector: 
assume that we have constructed the Hamiltonians $H(\s)$ corresponding to the different values 
of $N$ on the same probability space, and fix one realisation. From Proposition \eqv(prop14), 
we then know that, 
$$\limsup \frac 1 N \log \frac 1{\LL_{\pi_\b}(p)}\le \b^2(1+\z) \Eq(r10)$$
Let now $\l$ denote the spectral gap of $\EE$. $\l$ depends on the realisation of $H$ and on $N$. 
And let $\psi$ be the corresponding eigenvector. We assume that $\pi_\b(\psi^2)=1$. 
From [\rcite{FIKP}], we then know that 
$$\lim \frac 1 N \log \frac 1\l=\b\b_c \Eq(r11)$$
Therefore, provided we choose $\z$ small enough, we will have 
$$\limsup \frac 1 N\log \frac \l{\LL_{\pi_\b}(p)}\le -a \Eq(r12)$$, 
where $a>0$ is a deterministic constant that depends on $\b$. It then follows from \eqv(r4) 
that $\pi_\b(\vert \psi\vert) \le \exp(-a N)$ for large enough $N$ and with possibly a 
different value for the constant $a$. In other words the eigenvector $\psi$ becomes 
concentrated on its support. As a matter of fact, this is only another way to understand 
the fact that thermalisation times depend a lot on the initial law: eigenvectors corresponding 
to low eigenvalues become singular.

\vskip.5truecm

{\bf Proof of Theorem:\eqv(theo5)} 
 recalling \eqv(2.51), \eqv(4.16), \eqv(4.22) and \eqv(4.471), we get
$$
\eqalign{
\frac 1N \log T_N(\e,c,\eta) &\le
\frac 1N \log C_p+\frac{2-p}{pN}\log\frac 1\e+
\frac {4\log N}{N} \cr
& \quad +2\b\b_c 
\left( c_1(1+c) \frac {\log N}{N}\right)^{1/2} + \b^2 (1+\z) \cr
}\Eq(4.61)
$$ 
where $C_p$ is the constant in \eqv(2.51). (Remember that $d=\b(1+\z)$). 
Now taking first the limit $N\uparrow \infty$, we get 
$$\limsup \frac 1N \log T_N(\e,c,\eta) \le \b^2(1+\z)\Eq(4.61.1)$$
\eqv(4.61.1) is satisfied for 
all $\z>0$. (Just choose $p$ small enough so that \eqv(4.15) is satisfied). 
Therefore 
$$\limsup \frac 1N \log T_N(\e,c,\eta) \le \b^2$$
\qed

{\bf Proof of Theorem:\eqv(theo6)} 

The proof is a little more involved than the previous one.  
Choose 
$$\log\frac 1\e=N^{1/4}(\log N)^{3/4}$$
$$\z^2=12 \frac {\b_c}{\b} (c_1(1+c)\frac{\log N}N)^{1/2}$$
$\r=3/4$ and 
$$\frac p{1-p}=\frac 23\frac{\b^2}{\b_c^2+\b^2} \z^2$$
Then \eqv(4.15) and \eqv(4.151) are satisfied. Also 
$$\frac 2p=\frac 14 
\frac{\b^2+\b_c^2}{\b\b_c} \frac 1{\sqrt{c_1(1+c)}} \frac 1{\sqrt{N\log N}}$$
and we deduce the upper 
bound \eqv(4.901) from \eqv(4.61). 
The proof of \eqv(4.902) is similar, with now 
$\log(1/\e)=\d\log N$.\qed

\vskip1cm
\chap{V. The Medium  from the point of view of the process}5
\vskip.5truecm
\numsec=5
\numfor=1
\numtheo=1

\def\ib{\underline{i}}
In this section, we shall consider the process of the environment as seen from the
particle. This process will be denoted by $\o_t$.
For any fixed $N$, let $S_N\equiv \{-1,+1\}^N$. We endow $S_N$ with its natural group
structure i.e. for $\s,\s'\in S_N$, we let $\s.\s'\in S_N$ be the configuration
$(\s.\s')_i=\s_i\s'_i$. Let $\1$ be the configuration $(\1)_i=1$ for all $i$. For $1\le
i\le N$, we also define $\ib$ to be the configuration whose $i$-th coordinate is $-1$,
and the other coordinates are $+1$. Thus $\s.\ib$ is the configuration obtained by
flipping the $i$-th coordinate of $\s$.

Without loss of generality, we may, and will assume that our random  Hamiltonian $H$ is
defined on the canonical space $\O\equiv\R^{S_N}$. $\QQ$ is therefore the centered
product Gaussian probability on $\O$ of variance $N$. By duality, $S_N$ acts on $\O$
through the rule $(\s.h)(\s')\equiv h(\s.\s')$, where $\s,\s'\in S_N$ and $h\in \O$. 

For each choice of $H\in \O$, let us denote by $X^H$ the Metropolis dynamics with
Hamiltonian $H$, i.e. $X^H$ is the Markov process with generator 

$$L^H f (\s)\equiv 
\frac 1N \sum_{i=1}^N \e^{-\b [H(\ib.\s)-H(\s)]^+} 
\left( f(\ib.\s)-f(\s)\right)
\Eq(7.2)
$$ 
here, as before $[x]^+$ is the positive part of $x$.
We denote by  $P^H_t\equiv e^{t L^H}$, its semi-group, and  let $\E^H_\s$ be the 
law of $X^H$ when $X^H(0)=\s$.

Let us now define the stochastic process $\o_t\equiv X^H_t.H$. The state space of
$\o_t$ is $\O$. $\o_t$ is simply the Hamiltonian translated according to the position
of the particle. For instance note that, by definition, 
$\o_t(\1)=X^H_t.H(\1)=H(X^H_t)$ is nothing but the value of the Hamiltonian evaluated
at the position of the particle at time $t$. 
We consider the canonical construction of the Markov process, $X_t^H$, so we call 
$X_t$ the coordinate process on the space of cad-lag functions taking value in $S_N$.
We call $\E_\s^H$, the law of the Markov process with generator $L^H$ starting from $\s$.
We denote $e_\s^H$, the law of the process $\o_t\equiv X_t.H$ when $X_t$ is distributed according to $\E^H_\s$. 

By definition we have 
$$e^H_\s[\phi_1(\o_{t_1})...\phi_k(\o_{t_k})]
=\E^H_\s[\phi_1(X_{t_1}.H)...\phi_k(X_{t_k}.H)]
\Eq(7.20001)
$$

 The point is the following 

\Lemma (lem70)
$$
e^H_\s=e^{\s.H}_{1}
\Eq(7.200)
$$

\proof

Let $\phi$ be some measurable function on $\O$. Define $\phi^H(\s)\equiv \phi(\s.H)$.
Note that 
$$L^{\s.H}\phi^{\s.H}(\s')=L^H\phi^H(\s.\s')$$
this follows from
$$
\eqalign{
\left( L_N^{\s.H}f^{\s.H}\right)(\s')&=
\sum_{i=1}^N e^{-\b   [ \s.H(\underline i.\s')-\s.H(\s')]^+}
\left ( \phi\big((\underline i.\s').\s.H\big)-\phi\big(\s'.\s.H\big)\right)
\cr
&= 
\sum_{i=1}^N e^{-\b   [  H(\underline i.(\s.\s')- H(\s.\s')]^+}
\left ( \phi\big(\underline i.(\s  .\s').H\big)-\phi\big((\s .\s').H\big)\right)
\cr
&= 
\left( L_N^{ H}\phi^{ H}\right)( \s.\s')\cr
}\Eq(7.11)
$$
 
 Therefore, since $ \phi^{\s.H}(\s')=\phi^{H}(\s.\s')$, we get
$$
\left(e^{tL_N^{\s.H}}\phi^{ \s.H}\right)( \s')=\left( e^{tL_N^{ H}}\phi^{ H}\right)( \s.\s')
\Eq(7.12)
$$
Applying this last equality for $\s'=\1$, we have proved that
$$
e_{1}^{\s.H}\left[ \phi(\o_t)\right]=e_\s^H\left[ \phi(\o_t)\right]
\Eq(7.120)
$$
that is \eqv(7.200) holds for functions of one coordinate.

To extend it to an arbitrary cylindrical function we have, assuming $t_1<t_2$
$$
\eqalign{
e_{\s}^H\left[ \phi_1(\o_{t_1})\phi_2(\o_{t_2})\right]&
=\E^H_\s\left[\phi_1(X_{t_1}.H) \phi_2(X_{t_2}.H)\right]\cr
&=\E^H_\s\left[\phi_1(X_{t_1}.H) 
\E_{X_{t_1}}\left[\phi_2(X_{t_2-t_1}.H)\right]\right]\cr
}\Eq(7.121)
$$
where at the last step we have used that $X_t\equiv X_t^H$ is an homogeneous Markov process.
Using \eqv(7.200), we have
$$
\E_{X_{t_1}}^H\left[\phi_2(X_{t_2-t_1}.H)\right]=
\E_{1}^{X_{t_1}.H}\left[\phi_2(X_{t_2-t_1}.X_{t_1}.H)\right]
\Eq(7.122)
$$
using again \eqv(7.200) twice, we have also
$$
\eqalign{
&\E^H_\s\left[\phi_1(X_{t_1}.H) 
\E_{1}^{X_{t_1}.H}\left[\phi_2(X_{t_2-t_1}.X_{t_1}.H)\right]\right]\cr
&=\E^{\s.H}_1\left[\phi_1(X_{t_1}.\s.H) 
\E_{1}^{X_{t_1}.\s.H}\left[\phi_2(X_{t_2-t_1}.X_{t_1}.\s.H)\right]\right]\cr
&=\E^{\s.H}_1\left[\phi_1(X_{t_1}.\s.H) 
\E_{X_{t_1}}^{\s.H}\left[\phi_2(X_{t_2-t_1} .\s.H)\right]\right]\cr
}\Eq(7.123)
$$

Using once again the Markov property for $X_t$, we get
$$
\eqalign{
&\E^{\s.H}_1\left[\phi_1(X_{t_1}.\s.H) 
\E_{X_{t_1}}^{\s.H}\left[\phi_2(X_{t_2-t_1} .\s.H)\right]\right]\cr
&=\E^{\s.H}_1\left[\phi_1(X_{t_1}.\s.H) 
 \phi_2(X_{t_2 }.\s.H)\right] \cr
&=e_1^{\s.H}\left[\phi_1(\o_{t_1})\phi_2(\o_{t_2})\right]\cr 
}\Eq(7.124)
$$
Now it is easy to generalize what we just did to an arbitrary product 
of functions of one coordinate, then to cylindrical function and to measurable
function by the monotone class theorem. This ends the proof of the lemma. \qed

Note that $\o_t$ is the image of $X_t$ by the map $X_t \rightarrow X_t.H$.
In general the image of a Markov process is not Markovian, however here we have the

\Lemma(lem701){\it $\o_t$ is an homogeneous Markov process.}

\proof
It is enough to prove that
$$
 e^{\s.H}_1\left[\phi_1(\o_{t_1})\phi_2(\o_{t_2})\right]
  =
e_1^{\s.H}\left[\phi_1(\o_{t_1})
e^{\o_{t_1}}_1\left[\phi_2(\o_{t_2-t_1})\right]\right]
 \Eq(7.125)
$$
We have
$$
e^{\s.H}_1\left[\phi_1(\o_{t_1})\phi_2(\o_{t_2})\right]=
\E_{\s}^H\left[\phi_1(X_{t_1}.H) \phi_2(X_{t_2}.H)\right]
\Eq(7.126)
$$
Since $X_t$ is an homogeneous Markov process, we have
$$
\P^H_\s\left[ X_{t_1}=\s_1, X_{t_2}=\s_2\right]=
\P^H_\s\left[X_{t_1}=\s_1\right] 
\P^H_{\s_1}\left[X_{t_2-t_1}=\s_2\right]
\Eq(7.127)
$$
Therefore we get
$$
\eqalign{
&\E_{\s}^H\left[\phi_1(X_{t_1}.H) \phi_2(X_{t_2}.H)\right]=\cr
&\quad\sum_{H'_1} \phi_1(H'_1)\sum_{\s_1} \1_{\{\s_1.H=H'_1\}}
\P_{\s}^H\left[X_{t_1}=\s_1\right]\E^H_{\s_1}\left[\phi_2(X_{t_2-t_1}.H)\right]\cr
}\Eq(7.128)
$$ 
The point is that using \eqv(7.200), we have
$$
\eqalign{
&\sum_{\s_1} \1_{\{\s_1.H=H'_1\}}
\P_{\s}^H\left[X_{t_1}=\s_1\right]\E^H_{\s_1}\left[\phi_2(X_{t_2-t_1}.H)\right]\
=\cr
&\quad\quad e_1^{H'_1}\left[  \phi_2(\o_{t_2-t_1})\right] 
\sum_{\s_1} \1_{\{\s_1.H=H'_1\}}\P_{\s}^H\left[X_{t_1}=\s_1\right]\cr
}\Eq(7.129)
$$
Therefore we get
$$
\eqalign{
&\E_{\s}^H\left[\phi_1(X_{t_1}.H) \phi_2(X_{t_2}.H)\right] \cr
&\quad = 
\sum_{H'_1} \phi_1(H'_1)
e_1^{H'_1}\left[  \phi_2(\o_{t_2-t_1})\right] 
\sum_{\s_1} \1_{\{\s_1.H=H'_1\}}\P_{\s}^H\left[X_{t_1}=\s_1\right]\cr
&\quad=
\E_{\s}^H\left[\phi_1(X_{t_1}.H) e_1^{X_{t_1}.H}
\left[\phi_2(\o_{t_2-t_1})\right]\right]\cr
&\quad= e_1^{\s.H}\left[\phi_1(\o_{t_1})e^{\o_{t_1}}_1\left[\phi_2(\o_{t_2-t_1})\right]
\right]\cr
}\Eq(7.130)
$$
which is what we wanted to prove.\qed

Let $\p_\b^H$ be the Gibbs measure with Hamiltonian $H$ i.e. 
$$\p_\b^H(\s)\equiv \frac {e^{-\b H(\s)}}{Z^H(\b)}$$
and let us define the probability $\nu_\b^H$ on $\O$ by
$$
\nu_\b^H(f)\equiv 
\sum_{\s\in S_N} f(\s.H)\pi_\b^H(\s)
\Eq(7.6)
$$
when $f:\O\rightarrow \R$. That is for all $H'\in \O$
$$
\nu_\b^H(H')=\sum_{\s\in S_N} \pi_\b^H(\s) \1_{\{H'=\s.H\}}
\Eq(7.61)
$$

 We have the 

\Lemma(lem72){\it For each $H\in \O$, $ \nu_\b^H$ is an invariant 
and reversible measure for $\o_t$}

\proof

The invariance follows from
$$
 \eqalign{ 
\sum_{H'}\nu_\b^H(H') e_1^{H'}\left[\phi(\o_t)\right]=
&\sum_{H'}\sum_{\s\in\S_N} \pi_\b^H(\s) 
\1_{\{H'=\s.H\}}e_1^{H'}\left[\phi(\o_t)\right]\cr
&=\sum_{\s\in S_N} \pi_\b^H(\s)\E_{\s}^{H}\left[\phi^H(X_t)\right] \cr
&=\sum_{\s\in S_N}  \pi_\b^H(\s) \phi^H(\s)\cr
&=\nu^H_\b(\phi)\cr
}\Eq(7.62)
$$
where we have used the fact that $\pi_\b^H$ is invariant for $X_t$ at the third step.

The reversibility follows from
$$
 \sum_{\s\in S_N} \phi^H(\s) \pi_\b^H(\s) \E_\s^H\left[\psi^H(X_t)\right] 
  = 
\sum_{\s\in S_N} \psi^H(\s) \pi_\b^H(\s) \E_\s^H\left[\phi^H(X_t)\right]
 \Eq(7.63)
$$
 since  $\pi_\b^H$ is reversible for $X_t$. This ends the proof of the lemma.\qed

Now, for any bounded measurable function $f$ defined on $\O$, we have, as $t$ tends
to $+\infty$, 
$$ e^H_\s[f(\o_t)]=\E^H_\s[f(X^H_t.H)]\ra\nu_\b^H(f)\Eq(7.5)$$
We are interested in estimating the speed on convergence in \eqv(7.5). 
A fundamental fact is stated in the following lemma

\noindent {\bf \Lemma(lem71)} {\it For any $\varphi:\O_N \rightarrow \R$, 
$ \QQ\left[  e^H_{\s}(\varphi(\o_t))\right]$ is independent of $\s \in S_N$}
\medskip
\proof
This follows from the fact that on the one hand, for all $ \s \in S_N$ 
and for all $f: \O \rightarrow \R$, we have
$$
\QQ\left[f(H)\right]=\QQ\left[f( \s H)\right]
\Eq(7.9)
$$
since $\QQ$ is invariant by any permutation  of the configurations $H$.

Therefore, using \eqv(7.11), we have, for all $\varphi:\O  \rightarrow \R$
$$
\QQ\left[e^H_{\s}(\phi)\right]=\QQ\left[e_1^{\s.H}(\phi)\right]
=\QQ\left[e^H_1(\phi)\right]
\Eq(7.100000)
$$
which is what we wanted to prove.\qed

Now we can define the following time:
$$
T_{av}(\e)\equiv \inf \left\{
t>0\ \,  s.t.\  \sup_{s\geq t} \sup_{\varphi:\Vert\varphi\Vert_\infty\le 1}
\QQ\left[\left| 
e^H_\s(\varphi(\o_t))- \nu^{H}_{\b}(\varphi)\right| \right]
\le \e
\right\}
\Eq(7.15)
$$
here $\Vert \varphi\Vert_\infty=\sup_{\o\in \O}\vert \varphi(\o)\vert$.
$ T_{av}(\e)$ is the time such that  the average over the medium  
of the medium as seen 
from the process  is definitively within $\e$ of  the reversible 
measure $ \nu^{H}_{\b}$.
The main result of this section is

{\bf \Theorem(theo72)} {\it For all $\e>0$, for all $\b\le \b_c$,
$$
\limsup_{N \uparrow \infty} \frac 1N \log T_{av}(\e) \le \b^2
\Eq(7.16)
$$
}
{\proof}
Using Lemma \eqv(lem71), 
denoting by $d\eta(x)$ the uniform measure on $S_N$, we get 
$$
\QQ\left[\left| e^H_\s(\varphi(\o_t))- 
\nu^{H}_{\b}(\varphi)\right| \right]
= \int d\eta(\s)
\QQ\left[\left|  e^H_\s(\varphi(\o_t))- 
\nu^{H}_{\b}(\varphi)\right| \right]
\Eq(7.17)
$$
since the left hand side does not depends on $\s$.
Now using Tonelli's theorem we get
$$
 \int d\eta(\s)
\QQ\left[\left|   e^H_\s(\varphi(\o_t))- \nu^{H}_{\b}(\varphi)\right| \right]
=
\QQ\left[ 
\int d\eta(\s)\left|   e^H_\s(\varphi(\o_t))- 
\nu^{H}_{\b}(\varphi)\right| \right]
\Eq(7.18)
$$
 Now, since 
$$
  e^H_\s(\varphi(\o_t))=\E^H_\s(\varphi^H(X_t))
\Eq(7.19)
$$
using \eqv(7.6), we get 
$$
 \nu^{H}_{\b}(\varphi)= \pi^{H}_{\b}(\varphi^H)
\Eq(7.20)
$$
 Therefore
$$
\eqalign{
  e^H_\s(\varphi(\o_t))- \nu^{H}_{\b}(\varphi)&=
\E^H_\s(\varphi^H(X_t))- \pi^{H}_{\b}(\varphi^H)\cr
&= \left(P^{H}_{t}(\varphi^{H})\right)(\s)-  \pi^{H}_{\b}(\varphi^H)\cr
}
\Eq(7.21)
$$
therefore collecting what we just did  we get
$$
\eqalign{
\QQ\left[\left|   e^H_\s(\varphi(\o_t))- 
\nu^{H}_{\b}(\varphi)\right| \right]
&=
\QQ\left[\int d\eta(\s')\left| \left(P^{H}_{t}(\varphi^{H})\right)(\s')-  
\pi^{H}_{\b}(\varphi^H)\right|\right]\cr
}\Eq(7.22)
$$

To continue, recalling Proposition \eqv(prop11) and \eqv(prop13), for all $c$,
let $\AA(c)$ be the subspace of $\O$, of $\QQ$-probability bigger
than $1-2e^{-cN}$, which is  the intersection  of the two subspaces where we 
have the estimates \eqv(4.22) and \eqv(4.471). 

 Then we get 
$$
\eqalign{
&\QQ\left[\int d\eta(x)\left| \left(P^{H}_{t}(\varphi^{H})\right)(x)
 -\pi^{H}_{\b}(\varphi^H)\right|\right]
\le 2 \QQ\left[ \Vert \varphi \Vert_\infty \1_{\AA^c(c)}\right] +\cr
&\ \ \ \ \ + C_p t^{-p/(2-p)}
\QQ\left[ \1_{\AA(c)}
\Vert \varphi^H\Vert_\infty
(\LL^H_\h(p))^{-p/2} (\LL^H(p))^{-p^2/(4-2p)}\right]\cr
}\Eq(7.23)
$$
where the first part of the inequality follows 
from the fact that for all $H\in \O$ and all
$t>0$,  $P^{H}_t$ is a contraction operator 
from $L^{\infty}[\O,\R]$ into itself
and $\pi_{\b}^H$ is a probability measure. 
The second part follows from \eqv(2.6).
We  recall that $C_p=e^{-p/2} 2^{p/2} ((2-p)/p)^{-p^2/(4-2p)}$. 

 Using now Proposition \eqv(prop11), Proposition \eqv(prop13) and 
$\Vert \varphi \Vert_\infty\le 1$, we get
$$
\eqalign{
&\QQ\left[\left|   e^H_\s
(\varphi(\o_t))- \nu^{H}_{\b}(\varphi)\right| \right]
\le
2e^{-c N} +\cr
&+ C_p t^{-p/(2-p)}
(22)^{p^2/(4-2p)}
(4 )^{p/2} 
(N)^{8p/4-2p}( e^{2\b\b_c \sqrt{c_1(1+c) N\log N}}e^{\b dN})^{p ((4-p)/4-2p)}
}\Eq(7.24)
$$
 From now on the proof is exactly the same as the proof of Theorem \eqv(theo5).

At this point it is clear that we can also gives estimates that are similar to the
 ones given in Theorem \eqv(theo6) by using the same arguments as before
and the computation done in the proof of Theorem \eqv(theo6). 
Let us state it as a Theorem.

{\bf \Theorem(theo73)} {\it For all $N$ large enough, for all $\b\le \b_c$,
There exists a constant $c_1>0$, such that for all $c>0$, 
there exists a constant 
$C_0=C_0(c,\b)$ such that 
$$\eqalign{
&\frac 1N \log T_{av}(e^{- N^{1/4}(\log N)^{3/4}})\cr
\le &\b^2 +  
2\b \b_c \left({ \frac{c_1(1+c)\log N}{N}}\right)^{1/2} + c_2(\b,c)
\left( \frac {\log N}{N}\right)^{1/4}+C_0 (\frac{\log N}N)^{3/4}
}\Eq(7.25)
$$
where
$$
c_2(\b,c)\equiv \b\left(12 \b\b_c \sqrt{ c_1(1+c)}\right)^{1/2} 
+ \frac 14\frac{ \b^2+\b_c^2}{\b \b_c \sqrt{ c_1(1+c)}}
\Eq(7.26)
$$
Moreover  for all $\d>0$
$$
\eqalign{
&\frac 1N \log T_{av}(N^{-\d})
 \le \b^2 + 
\b\left( 12 \b\b_c\sqrt{c_1(1+c)}\right)^{1/2} (\frac{\log N}N)^{1/4}\cr
& +2\b \b_c\left({ \frac{c_1(1+c)\log N}{N}}\right)^{1/2} 
+\frac 14
\frac{\b^2+\b_c^2}{\b\b_c}\frac\d{\sqrt{c_1(1+c)}}(\frac{\log N}N)^{1/2}
+C_0 \frac{\log N}N
}\Eq(7.27)
$$
}


\vskip1cm

\chap{VI. Statics estimates for the  REM}6
\vskip.5truecm
\numsec=6
\numfor=1
\numtheo=1

 In this section we will give some estimates for the various
 constrained partition functions and partition functions on small
 spaces for the REM. These are just adaptations of similar estimates
 done in [\rcite{FIKP}] section 4.2.1.

 Let us first prove  Lemma {\eqv(lem12)}.
 We denote by $Z_\a(\b,\ge -d)\equiv Z_{j-1}(\b,\ge -d)[z_j,z_{>j}]$.
 Let $M$ be as in  Lemma {\eqv(lem12)}, and make the partition of the
 real interval $(-\infty, d N]$ with the intervals
$$
\D_0\equiv \left(-\infty, \b_c \frac {N}{M}\right]
\Eq(5.1)
$$
if $1\le k \le \frac d{\b_c} M-1$
$$
\D_k\equiv \left( \b_c \frac kM N,\b_c \frac {k+1}M \right]
\Eq(5.2)
$$
 Let 
$$
N_k=N_k(z_j,z_{>j})= \sum_{x\in \{-1,+1\}^{j-1}} \1_{\D_k}
( -H(x,z_j,z_{>j}))
\Eq(5.3)
$$
be the occupation number of the interval $\D_k$, it is easy to check
that, if $p_k=\P[-H(x) \in  \D_k]$, then 
$$
\b_c \frac {\sqrt N}M 2^{-\sfrac{(k+1)^2}{M^2} N}
< p_k <  
\b_c \frac {\sqrt N}M 2^{-\sfrac{k^2}{M^2} N}
\Eq(5.4)
$$
 Using the exponential Markov inequality and  optimizing we get
$$
\P\left[ N_k > \r_k \E(N_k)\right]
\le \exp\left\{ -\l_k 2^{\a N}\right\}
\Eq(5.5)
$$
where 
$$
\r_k=2^{N[ \sfrac{(k+1)^2}{M^2} -\a ]^+ +2}
\Eq(5.6)
$$
and if $\r_k p_k \ge 1$, $ \l_k=\infty$, while if $\r_k p_k<1$
$$
\l_k\equiv \r_k p_k \log \frac {\r_k(1-p_k)}{1-\r_k p_k}
-\log \left[ 1-p_k + \frac{\r_k p_k(1-p_k)}{1-\r_k p_k}\right]
\Eq(5.7)
$$

 It is not too long  to check that $\l_k \ge \r_k p_k c_1$ for some
positive constant $c_1$, and also $\r_k p_k \ge 2^{N/M^2}$, therefore
with our choice of $M$, we get
$$
\eqalign{
\P\left[ N_k > \r_k \E(N_k)\right]& \le 
\exp-\left\{ c_1{2^{N/M^2}}\right\}\cr
&\le 2^{-2N}\exp -(c N)
}\Eq(5.8)
$$
Note that the term $2^{-2N}$ will be more than enough to get uniformity
with respect to the index $i$ for the chosen family of path, the index $j$,
the configurations $z_j,z_{>j}$, and the index $k$.

Therefore, calling $A \equiv \sqrt{ \a M^2-1}$ and $D+1\equiv d
M/\b_c$ and using \eqv(5.4) and \eqv(5.6), we get 
$$ 
\eqalign{
Z_\a(\b,\ge -d) &\le 2^{j}e^{\b\b_c \sfrac NM}
+\sum_{k=1}^{A\wedge D}\frac {\sqrt{ N}}{M} 
e^{N(\a-\sfrac {k^2}{M^2})\sfrac{\b_c^2}{2}
+ \b\b_c \sfrac{(k+1)}{M} N} \cr
& \quad \quad 
+\sum_{k=A\wedge D +1}^D \frac {\sqrt {N}}{M}
e^{N\b\b_c \frac {(K+1)M}{ N }} 2^{N/M^2}\cr
}\Eq(5.9)
$$
where the last sum is not present if $D < A$.

 We have 
$$
2^{j}e^{\b\b_c \sfrac NM}\le
2^J e^{ \b\b_c \sqrt{ c_1(1+c)N\log N}}
\Eq(5.10)
$$

It is immediate to see that, if $D<A$ 
$$
\sum_{k=A\wedge D +1}^D \frac {\sqrt {N}}{M}
e^{N\b\b_c \frac {(K+1)M}{ N }} 2^{N/M^2}
\le 
c_u N^{3/2} e^{\b d N}
\Eq(5.11)
$$

 It remains to estimate the first sum in the right hand side of
\eqv(5.9).
Let us call it $ S(N)$, if we denote $x= K/M$, the maximum in the
exponent occurs for $x=\b/\b_c$.
Therefore, if $ A < \frac {\b}{\b_c} M$ we easily
get 
$$
\eqalign{
S(N) &\le 
\sqrt{N}e^{ \b\b_c \sqrt{ c_1(1+c)N\log N}} 
e^{N \b\b_c \sqrt
{\a}}\cr
& \le \sqrt{N}e^{ \b\b_c \sqrt{ c_1(1+c)N\log N}} e^{N \b d}\cr
}\Eq(5.12)
$$
where at the last step we have used $\b_c \sqrt{\a} < \b <d$

If $\frac {\b}{\b_c} M \le A < D $ we easily
get 
$$
S(N) \le \sqrt{N}e^{ \b\b_c \sqrt{ c_1(1+c)N\log N}} 
e^{N(\sfrac{\b^2}{2} +\a \sfrac{\b_c^2}{2})}
\Eq(5.13)
$$

If $ D \le A$, since $d>\b$, the maximum of the exponent occurs inside
the interval of summation therefore we easily
get that 
$$
S(N) \le \sqrt{N}e^{ \b\b_c \sqrt{ c_1(1+c)N\log N}} 
e^{N(\sfrac{\b^2}{2} +\a \sfrac{\b_c^2}{2})}
\Eq(5.14)
$$
collecting \eqv(5.10) to \eqv(5.14) we get \eqv(4.42) and \eqv(4.43).

The Lemma \eqv(lem10) is proved in exactly the same way, by making a
similar partition of $ [d N, +\infty)$, for proving \eqv(4.18).
Restricting the sum over $k$ to just the one corresponding to
$k= M \b/\b_c$, it is easy to get \eqv(4.19).


\vfill
\eject
\vskip1cm
\centerline{\bf References}
\vskip.3truecm

\item{[\rtag{B1}]}  Derrida B. (1980) 
{ Random Energy model: Limit of a family of
disordered systems.} {\it Phys. Rev. Lett} {\bf 45} 79--82.
\item{[\rtag{B2}]}  Derrida B. (1981) 
{ Random Energy model:
 An exactly solvable model of a spin glass.} 
{\it Phys. Rev. B} {\bf 24} 2613--2626. 

\item{[\rtag{Ei}]} Eisele, T.  (1983) 
{ On a Third Order Phase Transition.} {\it   Comm. Math. Phys.} {\bf 90}, 125--159. 
\item{[\rtag{FIKP}]} Fontes, L.R.G., Isopi, M, Kohayakawa, Y., and Picco, P. (1998)
{ The sprectral gap of the REM under metroplolis dynamics.}
{\it Ann. Appl. Prob} {\bf 8} 917--943.
\item{[\rtag{M}]} Mathieu P. (1997) 
{ Hitting times and spectral gap inequalities.} 
{\it Ann. Instit.  H.Poincar\'e. 
S\'erie Probabilit\'es} {\bf  33}. n 4,  437--465.
\item{[\rtag{M1}]} Mathieu P. (1999)   
{ Sur la convergence des marches al\'eatoires dans un milieu al\'eatoire et les
in\'egalit\'es de Poincar\'e g\'en\'eralis\'ees.} 
{\it CRAS Paris}, {\bf t. 329, s\'erie I}, 1015--1020 
\item{[\rtag{M2}]} Mathieu P. (2000) 
{ Convergence to equilibrium for spin glasses.} 
{\it Comm. Math. Phys. } {\bf 215}, 57--68. 
\item{[\rtag{OP}]} Olivieri E., Picco P. (1984) 
{ On the Existence of Thermodynamics for the Random Energy Model.} 
{\it  Comm. Math. Phys.} {\bf 96}, 125--144.
\item{[\rtag{GMP}]} Galves.A., Martinez S., and Picco P. (1989)
{ Fluctuations in Derrida's Random Energy and
 Generalized Random Energy Models.} 
{\it  J. Stat. Phys.} {\bf 54} 515--529.
\item{[\rtag{MPV}]} Mezard M., Parisi G., and Virasoro M.,(1980)
{\it Spin Glass theory and beyond}, { World Scientific}.
\item{[\rtag{SC}]} Saloff-Coste L. (1997).
{\it  Lectures on finite Markov chains.} 
Ecole d'\'et\'e de St-Flour. Lectures Notes in Math.
Springer Verlag, Berlin.

\bigskip

Pierre Mathieu \hfill\break
 CMI, Universit\'e de Provence,\hfill\break
39 Rue F. Joliot Curie, 13453 Marseille Cedex 13, France.\hfill\break
pmathieu$\!\,@$gyptis.univ-mrs.fr
\medskip
Pierre  Picco \hfill\break
 CPT.CNRS Luminy, case 907,\hfill\break 
13288 Marseille Cedex 9, France.\hfill\break 
Picco$\!\,@$cpt.univ-mrs.fr.

\end